\documentclass[11pt,a4paper,reqno]{amsart} 
\usepackage{amssymb,graphicx}

\usepackage{amssymb}
\def\theequation{\thesection.\arabic{equation}}

\newcommand{\bS}{{\bf \Psi}}
\newcommand{\bfA}{{\bf \Omega}}
\newcommand{\CP}{\mathbb{CP}}

\newcommand{\R}{\mathbb{R}}

\newcommand{\N}{\mathbb{N}}

\newcommand{\DG}{\mathcal{D}(\Gamma)}
\newcommand{\Rho}{\mathrm{P}}
\newcommand{\sol}{\mathcal{S}}
\newcommand{\koniec}{\begin{flushright}  $\Box $ \end{flushright}}
\newcommand{\intprod}{\;\rule{5pt}{.3pt}\rule{.3pt}{7pt}\;}

\def\p{\partial}
\def\be{\begin{equation}}
\def\ee{\end{equation}}
\newtheorem{theorem}{Theorem}[section]
\newtheorem{corollary}[theorem]{Corollary}
\newtheorem{lemma}[theorem]{Lemma}

\newtheorem{definition}[theorem]{Definition}
\begin{document}
\title{{\bf METRISABILITY OF TWO-DIMENSIONAL PROJECTIVE STRUCTURES}}
\author{Robert Bryant}
\address{
The Mathematical Sciences Research Institute\\
17 Gauss Way\\
Berkeley, CA 94720-5070\\USA.}
\email{bryant@msri.org}
\author{Maciej Dunajski}
\address{Department of Applied Mathematics and Theoretical Physics\\ 
University of Cambridge\\ Wilberforce Road, Cambridge CB3 0WA\\ UK.}
\email{m.dunajski@damtp.cam.ac.uk}
\author{Michael Eastwood}
\address{
School of Mathematical Sciences\\
Australian National University\\
ACT 0200, Australia.}
\email{meastwoo@member.ams.org}
\begin{abstract}{We carry out the programme of R. Liouville \cite{Liouville} to
construct an explicit local obstruction to the 
existence of  a Levi--Civita connection within a given projective 
structure $[\Gamma]$ on a surface. The obstruction is of order 5 in the 
components of a connection in a projective class. It can be expressed 
as a point invariant for a second order ODE whose integral curves are
the geodesics of $[\Gamma]$
or as a weighted 
scalar projective invariant of the projective class.  
If the obstruction vanishes we find the sufficient conditions for the 
existence of a metric in the real analytic case. 
In the generic case they are expressed  by 
the vanishing of two invariants of order
6 in the connection. 
In degenerate cases the sufficient obstruction is of order at most 8.}
 
\end{abstract}
\maketitle 

\section{Introduction}
Recall that a {\em projective structure} \cite{C22, Thomas, Eisenhart} 
on an open set
$U\subset \R^n$ is an equivalence class
of torsion free connections $[\Gamma]$. Two connections $\Gamma$ and
$\hat{\Gamma}$
are projectively equivalent if they share the same unparametrised geodesics.
This means that the geodesic flows project to the same foliation
of $\mathbb{P}(TU)$. The analytic expression for this equivalence class
is
\begin{equation}
\label{equivalence}
\hat{\Gamma}_{ab}^c= \Gamma_{ab}^c+\delta_a^c\omega_b+\delta_b^c\omega_a,
\qquad a, b, c=1, 2, ..., n
\end{equation}
for some one-form $\omega=\omega_a dx^a$.
A basic unsolved problem in projective differential geometry 
is to determine the  explicit criterion for
the {\em metrisability} of projective structure, i.e. answer the following 
question:
\begin{itemize}
\item 
What are the necessary and sufficient local conditions on a connection
$\Gamma_{ab}^c$ for the  existence of a one form  $\omega_a$ and
a symmetric non-degenerate tensor $g_{ab}$ such that the projectively
equivalent connection
\[
\Gamma_{ab}^c+\delta_a^c \omega_b+\delta_b^c\omega_a
\]
is the Levi-Civita connection for $g_{ab}$. (We are allowing 
Lorentzian metrics.)
\end{itemize}
We shall focus on  local metrisability, i.e. the pair $(g, \omega)$ with
$\det{(g)}$ nowhere vanishing
is required
to exist in a neighbourhood of a point $p\in U$. 
This problem leads to a vastly overdetermined system of
partial differential equations for $g$ and $\omega$. There are $n^2(n+1)/2$ 
components in
a  connection, and $(n+n(n+1)/2)$ components in a pair $(\omega, g)$.
One could therefore naively expect $n(n^2-3)/2$ conditions on $\Gamma$.

In this paper we shall carry out the algorithm laid out by R. Liouville 
\cite{Liouville} to
solve  this problem when $n=2$ and $U$ is a surface\footnote{Let us stress 
that the 
`solution' here means an explicit criterion, given by 
vanishing of a set of invariants, which can be 
verified on any representative of 
$[\Gamma]$.}.
In the  two-dimensional case the 
projective
structures are equivalent to second order ODEs which are cubic in
the first derivatives. To see it consider the geodesic equations
for $x^a(t)= (x(t), y(t))$ and eliminate the parameter $t$ between the
two equations
\[
\ddot{x}^c+\Gamma^c_{ab}\dot{x}^a\dot{x}^b=v\dot{x}^c.
\]
This yields the desired ODE for $y$ as a function of $x$
\begin{equation}
\label{ODE1}
\frac{d^2 y}{d x^2}=\Gamma^1_{22}\Big(\frac{d y}{d x}\Big)^3
+(2\Gamma^1_{12}-\Gamma^2_{22})\Big(\frac{d y}{d x}\Big)^2
+(\Gamma^1_{11}-2\Gamma^2_{12})\Big(\frac{d y}{d x}\Big)-
\Gamma^2_{11}.
\end{equation}
Conversely, any second order  ODEs cubic in
the first derivatives
\begin{equation}
\label{ODE2}
\frac{d^2 y}{d x^2}=A_3(x, y)\Big(\frac{d y}{d x}\Big)^3
+A_2(x, y)\Big(\frac{d y}{d x}\Big)^2
+A_1(x, y)\Big(\frac{d y}{d x}\Big)+A_0(x, y)
\end{equation}
gives rise to some projective structure as the independent
components of $\Gamma^c_{ab}$ can be read off from the $A$s up 
to the equivalence (\ref{equivalence}).      
The advantage of this formulation is that the projective ambiguity 
(\ref{equivalence})
has been removed from the problem as the combinations of the 
connection symbols in the ODE (\ref{ODE1}) are independent of the choice of 
the one form
$\omega$. There are $6$ components in $\Gamma^c_{ab}$  and $2$
components in  $\omega_a$, but only $4=6-2$ coefficients 
$A_\alpha(x, y), \alpha=0, ..., 3$.
The diffeomorphisms of $U$ can be used to further eliminate 
$2$ out of these $4$ functions (for example to make the equation 
(\ref{ODE2}) linear in the first derivatives) so one can say that
a general projective structure in two dimensions depends on two
arbitrary functions of two variables. We are looking
for invariant conditions, so  we shall not make use of this
diffeomorphism freedom. 

We shall state our first result.
Consider the 6 by 6 matrix   given in terms of its
row vectors
\begin{equation}\label{DEinvariant}
{\mathcal{M}}([\Gamma])=
({\bf V}, D_a{\bf V}, D_{(b}D_{a)}{\bf V} )
\end{equation}
which depends on the functions $A_\alpha$ and their derivatives up to 
order five. The vector field  ${\bf V}:U\rightarrow \R^6 $ is given by 
(\ref{formula_for_V}), the expressions  
$ D_a{\bf V}=\p_a {\bf V}-{\bf V}\bfA_a$ 
are computed using the right multiplication by 6 by 6 
matrices ${\bfA}_1, {\bfA}_2$ given 
by {(\ref{connection1})} and $\p_a=\p/\p x^a$.
We also make a recursive definition 
$D_aD_bD_c...D_d{\bf V}=\p_a(D_bD_c...D_d{\bf V})
-(D_bD_c...D_d{\bf V})\bfA_a$.
\begin{theorem}
\label{main_theorem}
If the projective structure $[\Gamma]$ is metrisable then
\begin{equation}
\label{DEinvariant1}
\det{({\mathcal M}([\Gamma]))}=0.
\end{equation}
\end{theorem}
There is an immediate corollary
\begin{corollary}
If the integral curves of a second order ODE 
\begin{equation}
\label{general_ODE}
\frac{d^2 y}{d x^2}=\Lambda\Big(x, y, \frac{d y}{d x}\Big),
\end{equation}
are geodesics of a Levi-Civita connection then 
$\Lambda$ is at most cubic
in $dy/dx$ and {\em(\ref{DEinvariant1})} holds.
\end{corollary}
The expression (\ref{DEinvariant1}) is written in a relatively compact form 
using $({\bf V}, {\bfA}_1, {\bfA}_2)$. 
All the algebraic manipulations which are required
in expanding the determinant 
have been done using  MAPLE code which can be 
obtained from us on request.

We shall prove  Theorem \ref{main_theorem} in three steps.
The first step, already taken by Liouville \cite{Liouville}, 
is to associate a linear
system of four PDEs for three unknown functions with each metrisable connection.
This will be done in the next Section. The second step will be prolonging this
linear system. This point was also understood by Liouville  
although he did not carry out the explicit computations.
Geometrically this will come down to constructing a connection
on a certain rank six real vector bundle over $U$.
The non-degenerate parallel sections of this bundle are in
one to one correspondence 
the metrics whose geodesics
are the geodesics of the given projective structure.
 In the generic case,
the bundle has no parallel sections and hence
the projective structure does not come from metric.  
In the real analytic case
the projective structure for which there is a single parallel section
depends on one arbitrary function of two variables, up to diffeomorphism. 
Finally we shall obtain (\ref{DEinvariant1}) as the integrability conditions
for the  existence of a parallel section
of this bundle. 
This will be done in Section \ref{sec_prolongation}.

In Section \ref{sec_sufficient} we shall present some sufficient
conditions for metrisability. 
All considerations here will be in the real analytic category.
The point is that even if  $[\Gamma]$ is locally metrisable 
around every point in $U$, the global metric on $U$ may not exist in the smooth
category even in the  simply-connected  case. Thus no set of local obstructions
can guarantee metrisability of the whole surface $U$.
\begin{theorem}
\label{theorem_sufficient5}
Let $[\Gamma]$ be a real analytic projective structure such that
rank $({\mathcal{M}}([\Gamma]))<6$ on $U$ and there exist $p\in U$ such that
rank $({\mathcal{M}}([\Gamma]))=5$ and $W_1W_3-W_2^2\neq 0$ at $p$, 
where $(W_1, W_2, ..., W_6)$ spans the kernel of 
${\mathcal{M}}([\Gamma])$. Then $[\Gamma]$ is metrisable in a 
sufficiently small
neighbourhood of $p$
if the rank of a $10$ by $6$ matrix with the rows
\[
({\bf V},  D_a{\bf V},  D_{(a} D_{b)}{\bf V},  D_{(a} D_b D_{c)}{\bf V})
\] 
is equal to $5$. 
Moreover this rank condition holds  if and only if two relative 
invariants $E_1, E_2$ of order 6 constructed
from the projective structure vanish.
\end{theorem}
We shall explain how to construct these two additional invariants
and show that the resulting set of conditions,
a single 5th order equation  (\ref{DEinvariant1}) and
two 6th order equations $E_1=E_2=0$ form an involutive system
whose general solution depends on three functions of two variables.
In the degenerate cases when rank$({\mathcal M}([\Gamma]))<5$
higher order obstructions will 
arise\footnote{We shall always assume that the rank of
$\mathcal{M}([\Gamma])$ is constant in a sufficiently small neighbourhood of 
some $p\in U$.}: one condition of order
8 in the rank 3 case and one condition of order 7 in the rank 4 case.
If rank $({\mathcal M}([\Gamma]))=2$ there is always a four parameter family
of metrics. If rank $({\mathcal M}([\Gamma]))<2$ then $[\Gamma]$ 
is projectively
flat in agreement with a theorem of Koenigs \cite{Koenigs}.
In general we have
\begin{theorem}
\label{theorem_rank8}
A real analytic projective structure $[\Gamma]$ is metrisable in 
a sufficiently small neighbourhood
of $p\in U$  if and only if the rank of a $21$ by $6$ matrix with the rows
\[
{\mathcal M}_{\mbox{{\em max}}}=({\bf V},  D_a{\bf V},  D_{(a} D_{b)}{\bf V}, 
D_{(a} D_b D_{c)}{\bf V},
D_{(a} D_b D_{c}D_{d)}{\bf V} ,   D_{(a} D_b D_{c}D_{d}D_{e)}{\bf V}
)
\]
is smaller than 6 and there exists a vector ${\bf W}$ in the kernel
of this matrix such that  $W_1W_3-W_2^2$ does not vanish at  $p$.
\end{theorem}
The signature of the metric underlying a projective structure can be Riemannian
or Lorentzian depending on the sign of $W_1W_3-W_2^2$. In the generic
case described by Theorem \ref{theorem_sufficient5}
this sign can be found by evaluating the polynomial (\ref{polynomial_adj}) 
of degree 10 in the entries of  ${\mathcal M}([\Gamma])$ at $p$.

In Section \ref{examples} we shall construct various examples
illustrating the necessity for the  genericity assumptions that 
we have made. 
In Section \ref{sec_twistor} we shall discuss the twistor approach 
to the  problem. In this approach a real analytic projective structure on $U$
corresponds to a complex surface $Z$ having a family of rational curves with
self-intersection number one. The metrisability condition and the associated linear
system are both deduced from the existence of a certain anti-canonical divisor
on $Z$.
In Section \ref{sec_tractor} we shall present an alternative tensorial 
expression for (\ref{DEinvariant1}) in terms of the curvature of the projective
connection and its covariant derivatives.  In particular we will shall
show that a section of the 14th power of the canonical bundle of $U$
\[
\det{({\mathcal M})}([\Gamma])\,(dx\wedge dy)^{\otimes 14}
\]
is a projective invariant.
The approach
will be that of tractor calculus \cite{Eastwood}.

In the derivation of the necessary condition (\ref{DEinvariant1})
we assume that the projective structure 
$[\Gamma]$ admits continuous fifth derivatives. The discussion of 
the sufficient conditions  and considerations in 
Section \ref{sec_twistor} require $[\Gamma]$ to be real analytic.
We relegate some long formulae to the Appendix.

We shall finish this introduction with a comment about the formalism 
used in the paper. The linear system governing the metrisability problem
and its prolongation are 
constructed in elementary way in Sections
\ref{sec_linearsystem}--\ref{sec_prolongation} and in tensorial
{\it tractor} formalism in Section \ref{sec_tractor}. The resulting 
obstructions are always given by invariant expressions.
The machinery of the Cartan connection could of course
be applied to do the calculations invariantly from the very beginning. 
This is in fact how some of the results have been obtained \cite{bryant_notes}.
The readers familiar with the Cartan approach
will realise that the rank six vector bundles used in our paper are
associated to the SL$(3, \R)$ principal bundle of Cartan. Such readers should
beware, however, that the connection $D_a$ that we naturally obtain on such a
vector bundle is not induced by the Cartan connection of the underlying
projective structure but is a minor modification thereof, as detailed for
example in~\cite{East_Mat}.
Various weighted invariants on $U$, like (\ref{DEinvariant1}), 
are pull-backs of functions from the total space of Cartan's bundle.
 
\vspace{2ex} {\bf Acknowledgements.} 
The first author is supported
by the National Science Foundation via grant DMS-0604195.
The second author is grateful to  Jenya Ferapontov, Rod Gover,
Vladimir Matveev and Paul Tod for helpful discussions. He also thanks  BIRS
in Banff and ESI in Vienna for hospitality where some of this research was
done. His work was partly supported by Royal Society and  London Mathematical
Society grants. The third author is supported by the Australian Research Council.
\section{Linear System}
\label{sec_linearsystem}
Let us assume that the projective structure $[\Gamma]$ is metrisable.
Therefore there exist a symmetric bi-linear form
\begin{equation}
\label{metric}
g=E(x, y)dx^2+2F(x, y) dxdy+G(x, y)dy^2
\end{equation}
such that the unparametrised geodesics of $g$ coincide
with the integral curves of $(\ref{ODE2})$. The
diffeomorphisms can be used to eliminate two arbitrary functions
from $g$ (for example to express  
$g$ in isothermal coordinates) but we shall not use this freedom. 

We want to determine
whether the four functions $(A_0, ..., A_3)$ arise from three functions
$(E, F, G)$ so one might expect only one condition on the $A$s. 
This heuristic numerology is wrong and we shall demonstrate in
Section \ref{sec_sufficient} that three conditions are needed
to establish sufficiency in the generic case\footnote{Additional 
conditions would arise if we demanded that there be more than one metric 
with the same unparametrised geodesics.
In our  approach  this situation corresponds 
to the existence of  two independent parallel sections
of the rank six bundle over $U$. The corresponding metrics
were, in the positive definite case,  found by J. Liouville 
(the more famous of the two Liouvilles) and characterised 
by Dini. They are of the form
(\ref{metric}) where $F=0, E=G=u(x)+v(y)$ up to diffeomorphism.
Roger Liouville whose steps we follow in this paper was a younger 
relative of  Joseph and attended his lectures
at the Ecole Polytechnique.}.

We choose a direct route and express the equation for non-parametrised
geodesics of $g$ in the form (\ref{ODE2}). Using the Levi-Civita
relation
\[
\Gamma_{ab}^c=\frac{1}{2}g^{cd}\Big(\frac{\p g_{ad}}{\p x^b}+ 
\frac{\p g_{bd}}{\p x^a}-\frac{\p g_{ab}}{\p x^d}\Big)
\]
and formulae (\ref{ODE1}), (\ref{ODE2}) yields the following
expressions
\begin{eqnarray}
\label{conditions}
A_0&=&\frac{1}{2}\frac{E\p_yE-2E\p_x F+F\p_xE}{EG-F^2},\nonumber\\
A_1&=&\frac{1}{2}\frac{3F\p_yE+G\p_xE-2F\p_xF-2E\p_xG}{EG-F^2},\nonumber\\
A_2&=&\frac{1}{2}\frac{2F\p_yF+2G\p_yE-3F\p_xG-E\p_yG}{EG-F^2},\nonumber\\
A_3&=&\frac{1}{2}\frac{2G\p_yF-G\p_xG-F\p_yG}{EG-F^2}.
\end{eqnarray}
This gives a first order nonlinear differential operator
\be
\label{1st_order_op}
\sigma^0:J^1(S^2(T^*U))\longrightarrow
J^0(\mbox{Pr}(U))
\ee
which carries the metric to its associated projective structure.
This operator is defined on the first jet space of symmetric two-forms 
as it depends on the metric and its derivatives. It takes
its values in the affine rank 4 bundle $\mbox{Pr}(U)$
of projective structures
whose associated vector bundle   $\Lambda^2(TU)\otimes S^3(T^*U)$
arises as a quotient in the exact sequence
\[
0\longrightarrow T^*U\longrightarrow TU\otimes S^2(T^*U)
\longrightarrow \Lambda^2(TU)\otimes S^3(T^*U)\longrightarrow 0.
\]
This is a more abstract way of defining the equivalence relation
(\ref{equivalence}). We will return to it in Section \ref{sec_sufficient}.
The operator $\sigma^0$  is homogeneous of degree 
zero so rescaling a metric by
a constant does not change the resulting projective structure.

Following Liouville \cite{Liouville}
we set
\[
E=\psi_1/\Delta^2, \quad F=\psi_2/\Delta^2, \quad
G=\psi_3/\Delta^2, \qquad \Delta =\psi_1\psi_3-{\psi_2}^2
\]
and substitute into (\ref{conditions}). This  yields an overdetermined system
of four linear first order PDEs for 
three functions $(\psi_1, \psi_2, \psi_3)$ and  proves the following 
\begin{lemma}[Liouville \cite{Liouville}]
\label{Liouville_lemma}
A projective structure $[\Gamma]$ corresponding
to the second order ODE {\em(\ref{ODE2})}
is metrisable on a neighbourhood of a point $p\in U$ iff there exists functions
$\psi_i(x, y), i=1,2,3$ defined on a neighbourhood of $p$ such that
\[
\psi_1\psi_3-{\psi_2}^2
\] 
does not vanish at $p$ and such that the equations
\begin{eqnarray}
\label{linear_system}
\frac{\p\psi_1}{\p x}&=&\frac{2}{3}A_1\psi_1-2A_0\psi_2,\nonumber\\
\frac{\p \psi_3}{\p y}&=&2A_3\psi_2-\frac{2}{3}A_2\psi_3,\nonumber\\
\frac{\p \psi_1}{\p y}+2\frac{\p \psi_2}{\p x}&=&\frac{4}{3}A_2\psi_1
-\frac{2}{3}A_1\psi_2-2A_0\psi_3,\nonumber\\
\frac{\p \psi_3}{\p x}+2\frac{\p \psi_2}{\p y}&=&
2A_3\psi_1
-\frac{4}{3}A_1\psi_3+\frac{2}{3}A_2\psi_2
\end{eqnarray}
hold on the domain of definition.
\end{lemma}
This linear system forms a basis of our discussion of the 
metrisability condition. It has recently been used in
\cite{BMM} to construct a list of metrics on a two-dimensional surface
admitting a two-dimensional group of projective transformations.
Its equivalent tensorial form, applicable
in higher dimensions, is presented for example in \cite{East_Mat}.
We shall use this form in Section \ref{sec_tractor}.

Here is a way to `remember' (\ref{linear_system}).
Introduce the symmetric projective connection $\nabla^{\Pi}$ with connection
symbols
\be
\label{thomas_symbols}
\Pi_{ab}^c=\Gamma_{ab}^c-\frac{1}{n+1}\Gamma_{da}^d\delta_b^c
-\frac{1}{n+1}\Gamma_{db}^d\delta_a^c
\ee
where in our case $n=2$. 
Formula (\ref{equivalence}) implies that the symbols 
$\Pi_{ab}^c$ do not 
depend on a choice of $\Gamma$ is a projective class. They are related
to the second order ODE (\ref{ODE2}) by
\[
\Pi_{11}^1=\frac{1}{3}A_1, \quad \Pi_{12}^1=\frac{1}{3}A_2, \quad \Pi_{22}^1=
A_3, \quad \Pi_{11}^2=-A_0, \quad \Pi_{21}^2=-\frac{1}{3}A_1, \quad
\Pi_{22}^2=-\frac{1}{3}A_2.
\]
The projective covariant derivative is defined on one-forms by
${\nabla^{\Pi}}_a\phi_b=\p_a\phi_b-\Pi_{ab}^c\phi_c$ with natural extension
to other
tensor bundles. The Liouville system (\ref{linear_system}) is then equivalent to
\be
\label{twistor_eq}
{\nabla^{\Pi}}_{(a}\sigma_{bc)}=0,
\ee
where the round brackets denote symmetrisation and $\sigma_{bc}$ is a 
rank 2 symmetric tensor with components 
$\sigma_{11}=\psi_1, \sigma_{12}=\psi_2, \sigma_{22}=\psi_3$.

\vspace{2ex}
We shall end this Section with a historical digression.
The solution to the metrisability problem has  been reduced to finding 
differential  relations between $(A_0, A_1, A_2, A_3)$ when
(\ref{conditions}), 
or equivalently (\ref{linear_system}),
holds. These relations are required to be 
diffeomorphism
invariant conditions, so we are searching for {\em invariants}
of the ODE (\ref{ODE2}) under the point transformations
\begin{equation}
\label{point}
(x, y)\longrightarrow (\bar{x}(x, y), \bar{y}(x, y)).
\end{equation}
The point invariants of 2nd order ODEs have been extensively
studied by the classical differential geometers in
late $19$th and early $20$th  century. The earliest reference
we are aware of is the work of Liouville \cite{Liouville87, Liouville}, 
who constructed point invariants of 2nd order ODEs cubic in the
first derivatives (it is easy to verify that 
the `cubic in the first derivative' condition is itself invariant under
(\ref{point})). The most complete work was produced by 
Tresse (who was a student of Sophus Lie) in his dissertation 
\cite{Tresse}. Tresse studied the general case
(\ref{general_ODE})
and classified all point invariants of a given differential order.
The first two invariants are of order four
\[
I_0=\Lambda_{1111}, \qquad I_1=
D_x^2\Lambda_{11}-4D_x\Lambda_{01}
-\Lambda_1 D_x\Lambda_{11}
+4\Lambda_1\Lambda_{01}-3\Lambda_0\Lambda_{11}+6\Lambda_{00},
\]
where 
\[
\Lambda_0=\frac{\p \Lambda}{\p y}, \qquad \Lambda_1=\frac{\p \Lambda}{\p y'}, 
\qquad
D_x=\frac{\p }{\p x}+y' \frac{\p }{\p y}
+\Lambda\frac{\p}{\p y'}.
\]
Strictly speaking these are only relative invariants as
they transform with a certain weight under (\ref{point}). 
Their vanishing is however invariant. Tresse
showed that if $I_0=0$, then $I_1$ is linear in $y'$. This is 
the case considered by Liouville.  To make contact with the work of Liouville
we note that 
$I_1=-6L_1-6L_2y'$ where the expressions
\begin{eqnarray}
\label{liouville_inv}
L_1&=&\frac{2}{3}\frac{\p^2 A_1}{\p x\p y}-\frac{1}{3}\frac{\p^2 A_2}{\p x^2}
-\frac{\p ^2 A_0}{\p y^2}+A_0\frac{\p A_2}{\p y}+  
A_2\frac{\p A_0}{\p y}\nonumber\\
&&- A_3\frac{\p A_0}{\p x}-2A_0\frac{\p A_3}{\p x}
-\frac{2}{3}A_1\frac{\p A_1}{\p y}
+\frac{1}{3} A_1\frac{\p A_2}{\p x},\nonumber\\
L_2&=&\frac{2}{3}\frac{\p^2 A_2}{\p x\p y}-\frac{1}{3}\frac{\p^2 A_1}{\p y^2}
-\frac{\p ^2 A_3}{\p x^2}-A_3\frac{\p A_1}{\p x}-  
A_1\frac{\p A_3}{\p x}\nonumber\\
&&+ A_0\frac{\p A_3}{\p y}+2A_3\frac{\p A_0}{\p y}
+\frac{2}{3}A_2\frac{\p A_2}{\p x}-\frac{1}{3} A_2\frac{\p A_1}{\p y}
\end{eqnarray}
were constructed by Liouville who has also proved that
\[
Y=(L_1 dx+L_2 dy)\otimes (dx\wedge dy)
\] is a projectively
invariant tensor.

The following result was known to both Tresse and Liouville
\begin{theorem}[ Liouville \cite{Liouville87}, Tresse \cite{Tresse}]
\label{thm_TL}
The 2nd order ODE {\em(\ref{general_ODE})}  is trivialisable by point 
transformation (i.e. equivalent to $y''=0$) iff $I_0=I_1=0$, or, 
equivalently, if $\Lambda$ is at most cubic in $y'$ and $Y=0$.
\end{theorem}
We note that the separate vanishing of $L_1$ or $L_2$ is not invariant.
If  both $L_1$ and $L_2$ vanish 
the projective structure is flat is the sense described 
in Section \ref{sec_tractor}.

%
%

\section{Prolongation and Consistency}
\label{sec_prolongation}
{\bf Proof of Theorem \ref{main_theorem}.}
The obstruction (\ref{DEinvariant1}) will arise as the 
compatibility condition for the system  (\ref{linear_system}).
This system 
is overdetermined, as there are more
equations than unknowns. We shall use the method of {\em prolongation}
and make (\ref{linear_system}) even more 
overdetermined\footnote{ 
Another approach more in the spirit of Liouville \cite{Liouville}
would be to eliminate $\psi_2$ and $\psi_3$ from (\ref{linear_system})
to obtain a system of two 3rd order PDEs for one function
$f:=\psi_1$
\[
(\p_x^3)f=F_1, \quad \p_y(\p_x^2)f=F_2,
\]
where $F_1, F_2$ are linear in $f$ and its first and second derivatives
with coefficients depending on $A_\alpha(x, y)$ and their derivatives
(the coefficient of $(\p_y^2)f$ in $F_1$ is zero).
The consistency $\p_y(\p_x)^3f=\p_x\p_y(\p_x)^2f$ gives a  linear 
equation for $\p_x(\p_y)^2f$. Then 
$\p_x(\p_y)^2\p_xf=(\p_y)^2(\p_x)^2f$ gives an equation for 
$(\p_y)^3f$. 
After this step the system is closed: all third order derivatives are
expressed in terms of lower order derivatives. To work out 
further consistencies impose 
$\p_x(\p_y)^3f=(\p_y)^3\p_xf$ which gives 
(when  all  3rd order equations are used) a second order linear PDE for $f$.
We carry on differentiating this second order relation to produce
the remaining
second order relations (because we  know all third order derivatives),
then the first order relations and finally an algebraic relation
which will constrain the initial data unless (\ref{DEinvariant1}) holds.} 
by 
specifying  the derivatives
of $\psi_i, i=1, 2, 3$ at any given point $(x, y, \psi_i)\in\R^5$,
thus determining a 
tangent plane to a solution surface (if one exists)
\[
(x, y)\longrightarrow (x, y, \psi_1(x, y), \psi_2(x, y), \psi_3(x, y)).
\]
For this we need six conditions, and  the system (\ref{linear_system}) 
consist of four equations. We need to add two conditions and we choose 
\begin{equation}
\label{munu}
\frac{\p \psi_2}{\p x}=\frac{1}{2}\mu, \qquad
\frac{\p \psi_2}{\p y}=\frac{1}{2}\nu,
\end{equation}
where $\mu, \nu$ depend on $(x, y)$.  
The 
integrability conditions $\p_x\p_y{\psi_i}=\p_y\p_x\psi_i$
give three PDEs for $(\mu, \nu)$ of the form
\begin{equation}
\frac{\p \mu}{\p x}=P, \qquad
\frac{\p\nu}{\p y}=Q, \qquad
\frac{\p \nu}{\p x} -\frac{\p\mu}{\p y}= 0,
\label{intermediat}
\end{equation}
where $(P, Q)$ given by (\ref{PandQ})
are expressions linear in $(\psi_i, \mu, \nu)$ with coefficients
depending on $A_\alpha$ and their $(x, y)$ derivatives. 

The system (\ref{intermediat}) 
is again overdetermined but
we still need to prolong it to specify the values
of all first derivatives. 
It is immediate that the complex 
characteristic variety of the system (\ref{linear_system}) is empty, so the
general 
theory (see Chapter 5 of \cite{Bryant_Chern}) implies that, after a finite
number of 
differentiations of these equations (i.e., prolongations), all of the 
partials of the $\psi_i$ above a certain order can be written in 
terms of lower order partials, i.e., the prolonged system will be 
complete. Alternatively, the Liouville system written in the form
(\ref{twistor_eq}) is one of the simplest examples covered by \cite{bceg} in
which the form of the prolongation is easily predicted.
In any case  
no appeal to the general theory is 
needed as it is easy to see that
completion is reached 
by adding one further equation
\begin{equation}
\label{Veq}
\frac{\p\mu}{\p y}=\rho,
\end{equation}
where $\rho=\rho(x, y)$ and imposing the consistency conditions on
the system of four PDEs  (\ref{intermediat}, \ref{Veq}). This leads
to 
\begin{equation}
\label{int_2}
\frac{\p\rho}{\p x}=R ,\qquad \frac{\p\rho}{\p y}= S,
\end{equation}
where $R, S$  given by (\ref{PandQ}) are functions of
$(\rho, \mu, \nu, \psi_i, x, y)$ which are linear in 
$(\rho, \mu, \nu, \psi_i)$.
After this step the prolongation process is finished and all the first
derivatives have been determined.
The final  compatibility
condition $\p_x\p_y \rho=\p_y\p_x \rho$ for the system (\ref{int_2}) 
yields
\begin{equation}
\label{compatibility}
\frac{\p R}{\p y}-\frac{\p S}{\p x}+S\frac{\p R}{\p \rho}
-R\frac{\p S}{\p \rho}=0.
\end{equation}
All the first derivatives are now determined, so  
(\ref{compatibility}) is an algebraic linear condition 
of the form
\begin{equation}
\label{linear_obstruction}
{\bf V}\cdot\bS:= \sum_{p=1}^6 V_p\Psi_p=0,
\end{equation}
where
\[
\bS=(\psi_1, \psi_2, \psi_3, \mu, \nu, \rho)^T
\]
is a vector in $\R^6$, and 
$
{\bf V}=(V_1, ..., V_6)  
$
where  
\begin{eqnarray}
\label{formula_for_V}
V_1&=& 2\frac{\p L_2}{\p y}+4A_2L_2+8A_3L_1,\quad
V_2= -2\frac{\p L_1}{\p y}-2\frac{\p L_2}{\p x}-\frac{4}{3}A_1L_2
+\frac{4}{3}A_2L_1,\nonumber\\
V_3&=&2\frac{\p L_1}{\p x}-8A_0L_2-4A_1L_1,\quad
V_4=-5L_2,\quad
V_5=-5L_1,\quad V_6=0
\end{eqnarray}
and $L_1, L_2$ are given by (\ref{liouville_inv})  .
We collect the linear PDEs 
(\ref{linear_system}, \ref{munu} \ref{intermediat}, \ref{Veq}, \ref{int_2})
as
\begin{equation}
\label{geometric_form}
d\bS+{\bfA}\,\bS=0,
\end{equation}
where
\[
{\bfA}={\bfA}_1\,dx+{\bfA}_2\,dy
\]
and $({\bfA}_1, {\bfA}_2)$ are $6$ by $6$ matrices with coefficients
depending on $A_\alpha$ and their first and second 
derivatives (\ref{connection1}). 
Now differentiate
(\ref{linear_obstruction}) twice with respect to $x^a=(x, y)$, and use
(\ref{geometric_form}). This yields six linear conditions
\begin{eqnarray}
\label{conditions11}
&&{\bf V}\cdot\bS=0,\\
&& (D_a {\bf V})\cdot\bS:= (\p_a {\bf V}-
{\bf V}\;{\bfA}_a)\cdot \bS=0,\nonumber\\
&&(D_b D_a {\bf V})\cdot\bS:=
(\p_b\p_a{\bf V}-(\p_b {\bf V})\;{\bfA}_a-(\p_a {\bf V})\;{\bfA}_b
-{\bf V}\;(\p_b {\bfA}_a-{\bfA}_a{\bfA}_b))\cdot\bS=0\nonumber
\end{eqnarray}
which must hold, or there are no solutions to (\ref{linear_system}). 
Therefore the determinant of the associated $6$ by  $6$ matrix  
(\ref{DEinvariant}) must vanish,
thus giving our first desired metrisability condition (\ref{DEinvariant1}).
We note that the expression $(D_b D_a {\bf V})\cdot\bS$ in
(\ref{conditions11}) is  symmetric in
its indices. This symmetry condition reduces
to ${\bf V} { F}=0$ (where ${F}$ is given by (\ref{six-curvature}))
and holds identically.

The expression $\det{({\mathcal M}([\Gamma]))}$ is $5$th order in the 
derivatives of connection coefficients. It does not vanish on a generic 
projective structure, but vanishes on metrisable connections (\ref{conditions})
by construction.
This ends the proof of Theorem \ref{main_theorem}.\koniec

In the next Section we shall need the following generalisation
of the symmetry properties of (\ref{conditions11}).
Let  $D_a{\bf W}=\p_a{\bf W}- {\bf W}{\bfA}_a$, where 
 ${\bf W}:U\rightarrow \R^6$. Then 
\[
[D_a, D_b]{\bf W}=({\bf W} F)\varepsilon_{ab}=W_6{\bf V}\varepsilon_{ab},
\]
where $\varepsilon_{00}=\varepsilon_{11}=0, 
\varepsilon_{01}=-\varepsilon_{10}=1$.
Thus 
\begin{eqnarray}
\label{symmetry_prop}
D_a D_b {\bf V}&=&D_{(a} D_{b)} {\bf V},\quad 
 D_a D_b D_c{\bf V}=D_{(a} D_b  D_{c)} {\bf V}+\varepsilon_{ab}L_c{\bf V},
\quad ...\,, \nonumber\\ D_{a_1} D_{a_2}... D_{a_k}{\bf V}&=&
D_{(a_1} D_{a_2}...  D_{a_k)} {\bf V}+o(k-2)
\end{eqnarray}
where $o(k-2)$ denotes terms linear in $D_{(a_1} D_{a_2} ... D_{a_m)}$
where $m\leq k-2$. Thus we can restrict ourselves to
the symmetrised expressions as the antisymmetrisations do not add any new
conditions.
\section{Sufficiency conditions} 
\label{sec_sufficient}
It is clear from the discussion in the preceding Section that
the condition (\ref{DEinvariant1}) is necessary for the existence of
a metric in a given projective class. It is however not sufficient and in
this Section we shall establish some sufficiency conditions in the real 
analytic case. We require the real analyticity in order to be able to apply the Cauchy--Kowalewski
Theorem to the prolonged system of PDEs. 
In particular Theorem \ref{cartan_test} which underlies our approach in this Section builds on the
Cauchy--Kowalewski Theorem.

Let us start off by rephrasing  the construction presented 
in the last Section
in the geometric language.
The  exterior differential ideal ${\mathcal I}$
associated to the prolonged system (\ref{geometric_form}) 
consists of six one-forms
\begin{equation}
\label{ideal_forms}
\theta_p=d\bS_p+(({\bfA}_a)_{pq}\bS_q)\, dx^a, \qquad
p, q=1, ..., 6\quad a=1,2.
\end{equation}
Two vector fields  
annihilating the one-forms $\theta_p$ span the solution
surface in $\R^8$.  The closure of this ideal comes down to
one compatibility  (\ref{linear_obstruction}).
We now want to find one parallel section $\bS:U\rightarrow \mathbb{E}$ 
of a rank six vector bundle $\mathbb{E}\rightarrow U$
with a connection $D= d+{\bfA}$. Locally the total space
of this bundle is an open set in $\R^8$.

Differentiating (\ref{geometric_form}) and eliminating $d\bS$ yields
${\bf F}\bS=0$, where
\begin{eqnarray*}
{\bf F}&=&d {\bfA}+{\bfA}\wedge {\bfA}=
(\p_x {\bfA}_2-\p_y {\bfA}_1+[{\bfA}_1, {\bfA}_2])dx \wedge dy\\
&=&Fdx\wedge dy
\end{eqnarray*}
is the curvature of $D$.
Thus we need
\begin{equation}
\label{holonomy1}
F\bS=0,
\end{equation}
where $F=F(x, y)$ is a $6$ by  $6$ matrix
given by (\ref{six-curvature}).
We find that this matrix is of rank
one and in the chosen basis its first five rows vanish and its bottom row
is given by the vector  ${\bf V}$ with components given by 
(\ref{formula_for_V}). Therefore (\ref{holonomy1}) 
is equivalent to (\ref{linear_obstruction}).
We differentiate the condition (\ref{holonomy1})
and use (\ref{geometric_form})  to produce algebraic matrix equations
\[
F\bS=0,\quad (D_aF)\bS=0,\quad (D_aD_bF)\bS=0,\quad (D_aD_bD_cF)\bS,\quad ...
\]
where $D_a F=\p_aF+[{\bfA}_a, F]$.
Using the symmetry argument 
(\ref{symmetry_prop}) shows that
after $K$ differentiations  this leads to $n(K)=1+2+3+...+(K+1)$
linear equations
which we write as
\[
{\mathcal F}_K \bS=0,
\]
where ${\mathcal F}_K$ is a $n(K)$ by $6$ matrix depending on $A$s
and their derivatives. We also set 
${\mathcal F}_0=F$.

We continue differentiating and adjoining 
the equations. 
The Frobenius Theorem  
adapted to  (\ref{holonomy1}) and (\ref{geometric_form})
tells us when we can stop the process.
\begin{theorem}
\label{cartan_test}
Assume that the ranks of the matrices  
${\mathcal F}_K, K=0, 1, 2, ...$ 
are maximal
and constant\footnote{This can always be achieved 
by restricting to a  sufficiently small neighbourhood of some point
$p\in U$.}.
Let $K_0$ be the smallest natural number such that
\begin{equation}
\label{rank_condition}
\mbox{rank}\; ({\mathcal F}_{K_0})=\mbox{rank}\; ({\mathcal F}_{K_0+1}).
\end{equation}
If $K_0$ exists then 
$\mbox{rank} ({\mathcal F}_{K_0})=\mbox{rank} ({\mathcal F}_{K_0+k})$ 
for $k\in\N$
and the space of parallel sections {\em(\ref{geometric_form})} of  $d+\bfA$
has 
dimension \[{\mathcal{S}}([\Gamma])=6-\mbox{rank} ({\mathcal F}_{K_0}).\]
\end{theorem}
The first and second derivatives of (\ref{holonomy1})
will produce six independent conditions on $\bS$, and these conditions 
are precisely
(\ref{conditions11}). Thus the necessary metrisability condition 
(\ref{DEinvariant1})  comes down to restricting the holonomy 
of the connection $D$ on the rank six vector bundle $\mathbb{E}$.

We shall now assume that (\ref{DEinvariant1}) holds and use 
Theorem \ref{cartan_test}
to construct the sufficient conditions for the existence of a Levi--Civita 
connection in a given projective class.
First of all there must exist a vector ${\bf W}=(W_1, ..., W_6)^T$ in the 
kernel of  ${\mathcal M}([\Gamma])$, such that 
$W_1  W_3-(W_2)^2\neq 0$. This will guarantee that the corresponding  
quadratic form (if one exists) on $U$ is non-degenerate.
It is straightforward to verify in the case
when ${\mathcal M}([\Gamma])$ has rank 5 as then 
$\mbox{kernel}\,({\mathcal M}([\Gamma]))$ is spanned by any non-zero column of
$\mbox{adj}\,({\mathcal M}([\Gamma]))$ where the adjoint of a matrix
$\mathcal{M}$ is defined by $\mathcal{M}\, \mbox{adj}(\mathcal M)=
\det{(\mathcal M)} \,I$.
The entries of $\mbox{adj}\,({\mathcal{M}}([\Gamma]))$ are determinants
of the co-factors of ${\mathcal{M}}([\Gamma])$ and thus
are polynomials
of degree 5 in the entries of ${\mathcal{M}}([\Gamma])$
so
\be
\label{polynomial_adj}
P([\Gamma])=W_1  W_3-(W_2)^2
\ee
is a polynomial of degree 10 in the entries of  ${\mathcal M}([\Gamma])$.
\begin{definition}
A projective structure for which {\em (\ref{DEinvariant1})} holds 
is called generic in a neighbourhood of $p\in U$
if rank ${\mathcal M}([\Gamma])$ is maximal and equal
to $5$ and
$P([\Gamma])\neq 0$ in this neighbourhood.
\end{definition}
In this generic case Theorem \ref{cartan_test} and Lemma \ref{Liouville_lemma}
imply that
there will exist a Levi--Civita connection 
in the projective class if the rank of the next
derived matrix ${\mathcal F}_3$
does not go up and is equal to five.
We shall see that this can be guaranteed by imposing
two more 6th order conditions on $[\Gamma]$. 

{\bf Proof of Theorem \ref{theorem_sufficient5}}. 
First note that, in the generic case, the three vectors
\[
 {\bf V},  \quad {\bf V}_a: = \p_a {\bf V}-{\bf V}\;{\bfA}_a, \qquad a=1,2 
\]
must be linearly independent or otherwise
the rank of ${\mathcal M}([\Gamma])$ would be at most 3. Now pick two 
independent vectors from the set
\[
{\bf V}_{ab}:= (\p_b\p_a{\bf V}-(\p_b {\bf V})\;
{\bfA}_a-(\p_a {\bf V})\;{\bfA}_b
-{\bf V}\;(\p_b {\bfA}_a-{\bfA}_a{\bfA}_b))
\]
such that the resulting set of five vectors is independent.
Say we have picked ${\bf V}_{00}$ and ${\bf V}_{11}$.
We now take the third derivatives of 
(\ref{linear_obstruction}) with respect to $x^a$ and use
(\ref{geometric_form}) to eliminate derivatives of ${\bf\Psi}$. 
This adds four vectors to
our set of five and so a priori we need to satisfy four six order
equations to ensure that the rank does not go up. However only
two of these are new and the other two are derivatives
of the 5th order condition  (\ref{DEinvariant1}). 
Before we shall prove this statement examining the images
of linear operators induced from (\ref{1st_order_op}) on jet spaces
let us indicate why this counting works. Let ${\bf V}_{ab...c}$
denote the vector in $\R^6$ annihilating 
$\bS$ (in the sense of (\ref{linear_obstruction})) which is obtained
by eliminating
the derivatives of $\bS$ from $\p_a\p_b...\p_c ({\bf V}\cdot\bS)=0$.
We have already argued in  (\ref{symmetry_prop}) that
the antisymmetrising over any pair of indices in 
${\bf V}_{ab...c}$ only adds lower order conditions. Thus we shall
always assume that these expressions are symmetric. We shall
also set ${\bf V}_0={\bf V}_x, {\bf V}_1={\bf V}_y$.

Our assumptions imply that 
\be
\label{my_dependence}
{\bf V}_{xy}= 
c_1{\bf V}+ c_2{\bf V}_x +c_3{\bf V}_y 
+c_4{\bf V}_{xx} +c_5{\bf V}_{yy}
\ee
for some functions $c_1, ..., c_5$ on $U$.
The two six order conditions
\be
\label{6th_orders}
E_1:=\det{
\left(
\begin{array}{c}
{\bf V}\\
{\bf V}_x\\
{\bf V}_y\\
{\bf V}_{xx}\\
{\bf V}_{yy}\\
{\bf V}_{xxx}
\end{array}
\right )},\quad
E_2:=\det{
\left(
\begin{array}{c}
{\bf V}\\
{\bf V}_x\\
{\bf V}_y\\
{\bf V}_{xx}\\
{\bf V}_{yy}\\
{\bf V}_{yyy}
\end{array}
\right )},\quad
\ee
have to be added for sufficiency. Now differentiating (\ref{my_dependence})
w.r.t $x, y$  and using  ${\bf V}_{xyy}={\bf V}_{yyx}, 
{\bf V}_{xyx}={\bf V}_{xxy}$ (which hold modulo
lower order terms), implies
that ${\bf V}_{xyy}$ and  ${\bf V}_{xyx}$
are in the span of 
$\{ {\bf V}, {\bf V}_x, {\bf V}_y ,{\bf V}_{xx}, {\bf V}_{yy}, {\bf V}_{xxx},
{\bf V}_{yyy} \}$
and no additional conditions need to be added. This procedure
can be  repeated if instead ${\bf V}_{yy}$ belongs to the span of
$\{ {\bf V}, {\bf V}_x, {\bf V}_y,{\bf V}_{xx},{\bf V}_{xy}\}$.
  
\vspace{2ex}
Now we shall present the general argument. Consider the homogeneous 
differential operator
(\ref{1st_order_op}). It maps the 1st jets of metrics on $U$
to the 0th jets of projective structures. Differentiating
the relations (\ref{conditions}) prolongs this operator to bundle maps
\be
\label{prolonged_operator}
\sigma^k:J^{k+1}(S^2(T^* U))\longrightarrow 
J^k(\mbox{Pr}(U))
\ee
from $(k+1)$-jets of metrics to $k$-jets of projective structures.
It has at least one dimensional fibre because of the  
homogeneity of $\sigma^0$. The rank of $\sigma^k$ is
not constant as we already know that the system (\ref{conditions}) 
(or its equivalent linear form (\ref{linear_system})) does not have to admit
any solutions in general but will admit at least one solution
if the projective structure is metrisable.
The table below gives the ranks of the jet bundles
of metrics and projective structures, the dimensions of the fibres
of $\sigma^{k}$ and finally the image codimension. 
The number of new conditions on $[\Gamma]$ arising at each step is denoted
by a bold
figure in the column $\mbox{co-rank}( \mbox{ker}\sigma^k)$.

\vspace{2ex}

\centerline{\small 
\noindent \begin{tabular}{|c|c|c|c|c|}  
\hline  $k$ & $\mbox{rank}(J^{k+1}(S^2(T^*U))) $ & 
$ \mbox{rank}(J^k(\mbox{Pr}(U)))$  &$\mbox{rank}( \mbox{ker}\sigma^k)$&
$\mbox{co-rank}( \mbox{ker}\sigma^k)$\\ 
\hline $-1$ & $3 $ & $ -$ & $-$ & $-$ \\
 \hline $0$ & $9 $ & $4 $ & $5$ & $0$ \\
\hline $1$ & $18 $ & $ 12$ & $6$ & $0$ \\
 \hline $2$ & $30 $ & $24 $ & $6$ & $0$ \\
 \hline $3 $ & $ 45 $ & $ 40 $ & $5 $ & $0 $ \\
\hline $4 $ & $ 63 $ & $60  $ & $3 $ & $0 $ \\
\hline $5 $ & $ 84 $ & $ 84 $ & $ 1$ & $1={\bf 1} $ \\
\hline $ 6$ & $ 108 $ & $112  $ & $1 $ & $5=3+{\bf 2} $ \\
\hline $ 7$ & $ 135 $ & $144  $ & $1 $ & $10=6+6-2 $\\
\hline
\end{tabular}}

\vspace{2ex}

There is no obstruction on a projective 
structure before the order 5
so $\sigma^k$  are onto and generically 
submersive for $k<4$.  At $k=5$ there has to be at least a
$1$-dimensional fiber, so the image of the derived map can at most be
$83$-dimensional at its smooth points.  
In fact, we have shown that there is a
condition there, given by (\ref{DEinvariant1}) , 
so it must define a codimension 1 variety
that is generically smooth. When the matrix ${\mathcal{M}}([\Gamma])$ 
has rank 5, the equation
  (\ref{DEinvariant1}) is regular, so it follows that, outside the region
where (\ref{DEinvariant1}) ceases to be a regular 5th order PDE 
the solutions of this PDE
will have their $k$-jets constrained by the derivatives of  
(\ref{DEinvariant1})  of order $k-5$
or less. This shows that, at $k=6$, the $6$-jets of the regular solutions of
(\ref{DEinvariant1})  will have codimension 3 in all 
6-jets of projective structures, i.e.,
they will have dimension $112-3 = 109$.  However, we know that the image of
the 7-jets of metric structures can have only dimension $108-1 = 107$. Thus,
the 6-jets of regular metric structures have codimension 2 in the 6-jets
of regular solutions of  (\ref{DEinvariant1}).  
That is why there have to be two more 6th
order equations
\be
\label{two_new_eq}
E_1=0, \qquad E_2=0.
\ee
The image in 6-jets has total
codimension 5, i.e., it is cut out by a 5th
order equation and four 6th order equations. However, two of the 6th order
equations are obviously the derivatives of the 5th order equation.  The
next line shows that, at 7th order, the image has only codimension 10, which
means that there must be 2 relations between the first derivatives of the
6th order equations and the second derivatives of the 5th order equation
which implies that the resulting system of three equations is involutive.
This ends the proof of Theorem \ref{theorem_sufficient5}.\koniec

The analysis of the non-generic cases where the rank of 
${\mathcal M}([\Gamma])<5$ is slightly more complicated.
The argument based on the dimensionality of jet bundles associated to
(\ref{prolonged_operator}) breaks down as the PDE
$\det{{\mathcal {M}}([\Gamma])}=0$ is not regular and does not define
a smooth co-dimension 1 variety in $J^5(\mbox{Pr}(U))$.

Let ${\mathcal{S}}([\Gamma])$ be the dimension of the vector
space of solutions to the linear system (\ref{linear_system}).
Some of these solutions may correspond to degenerate quadratic forms on $U$
but nevertheless we have
\begin{lemma}
If $\sol([\Gamma])>1$ then there are $\sol([\Gamma])$ independent 
non-degenerate quadratic forms among the solutions
to  {\em(\ref{linear_system})}.
\end{lemma}
{\bf Proof.} Let us assume that at least one solution
of (\ref{linear_system}) gives rise to a  quadratic form 
which is degenerate (rank 1) everywhere. We can choose coordinates
such that this solution is of the form $(\psi_1, 0, 0)$.
The statement of the Lemma will follow if we can show that
there is no other solution of the form 
$(\phi(x,y) \psi_1, 0, 0)$ where $\phi(x,y)$ is a non-constant
function. The Liouville system (\ref{linear_system}) is readily
solved in this case to give
\[
A_1=\frac{3}{2}\frac{1}{\psi_1}\frac{\p\psi_1}{\p x}, 
\quad A_2=\frac{3}{4}\frac{1}{\psi_1}\frac{\p\psi_1}{\p y}, \quad A_3(x, y)=0
\]
with $A_0$ unspecified. Thus for a given projective class
the only freedom in this solution is to rescale $\psi_1$ by a constant.
\koniec

{\bf Proof of Theorem \ref{theorem_rank8}.}
We shall list the number and the
order of obstructions one can expect depending on the rank
${\mathcal M}([\Gamma])$.
\begin{itemize}
\item If rank ${\mathcal M}([\Gamma])<2$ the projective structure
is projectively flat as $L_1=L_2=0$, and the second order
ODE is equivalent to $y''=0$ by Theorem \ref{thm_TL}. This is obvious if
rank ${\mathcal M}([\Gamma])=0$ as then ${\bf V}=0$ and
formula (\ref{formula_for_V}) 
gives $L_1=L_2=0$.

If rank ${\mathcal M}([\Gamma])=1$  then
\be
\label{gamma_1}
\p_a {\bf V}-{\bf V}{\bf \Omega}_a=\gamma_a{\bf V}
\ee
for some $\gamma_a$. Using  the expressions (\ref{connection1}) 
for ${\bf \Omega}_a$  and the formula (\ref{formula_for_V}) yields
\[
{\bf V}{\bf \Omega}_a=(*, *, *, *, *, 5L_a)
\]
where $*$ are some terms which need not concern us and $L_a$ are
the Liouville expressions (\ref{liouville_inv}). Combining this with
$(\ref{gamma_1})$ yields $L_1=L_2=0$.
\item
If rank ${\mathcal M}([\Gamma])=2$ then
\be
\label{rank_2}
{\bf V}+c_1{\bf V}_x+c_2{\bf V}_y=0
\ee
for functions $c_1, c_2$ at least one of which does not identically vanish.
Differentiating this relation and using the fact
that ${\bf V}_{ab}\in\mbox{span}\{ {\bf V}, {\bf V}_a\}$
we see that no new relations arise and so the system is closed
at this level. In this case there exists a four dimensional family of
metrics compatible with the given projective structure.
\item
If rank ${\mathcal M}([\Gamma])=3$ we have to consider two cases.
If $\{{\bf V}, {\bf V}_x, {\bf V}_y\}$ are linearly independent
then reasoning as above shows that
further differentiations do not add any new conditions.
The other possibility is that $\{{\bf V}, {\bf V}_x, {\bf V}_{xx}\}$
or $\{{\bf V}, {\bf V}_y, {\bf V}_{yy}\}$ are linearly independent.  
Let us concentrate on the first
case  (or swap  $x$ with $y$ if necessary). Taking further
$x$ derivatives may increase the rank of the resulting system, but
the $y$ derivatives will not yield any new conditions as can be seen
by mixing the partial derivatives and using 
\[
c_0{\bf V}+c_1{\bf V}_x+c_2{\bf V}_{y}=0,
\]
which is a consequence of the rank 3 condition.

Let us assume that
the rank increases to 5 by adding two vectors 
${\bf V}_{xxx}, {\bf V}_{xxxx}$ (otherwise the system is closed with rank
3 or 4). The rank will stay 5 if one further
differentiation does not add new conditions.
Thus the first and only obstruction in this case is
of order 8 in the projective structure
\be
\label{rank_8}
\det{
\left(
\begin{array}{c}
{\bf V}\\
{\bf V}_x\\
{\bf V}_{xx}\\
{\bf V}_{xxx}\\
{\bf V}_{xxxx}\\
{\bf V}_{xxxxx}
\end{array}
\right )}=0.
\ee
\item The analogous procedure can be carried over if 
rank(${\mathcal M}([\Gamma]))=4$. Assuming that the four linearly independent
vectors are  
$\{{\bf V}, {\bf V}_x, {\bf V}_{y},  {\bf V}_{xx}  \}$ leads 
to one obstruction of order 7 
 \[
\det{
\left(
\begin{array}{c}
{\bf V}\\
{\bf V}_x\\
{\bf V}_{xy}\\
{\bf V}_{xx}\\
{\bf V}_{xxx}\\
{\bf V}_{xxxx}
\end{array}
\right )}=0.
\]
\end{itemize}
This completes the proof of Theorem \ref{theorem_rank8}.\koniec

As a corollary from this analysis we deduce the result of Koenigs 
\cite{Koenigs}
\begin{theorem}\cite{Koenigs}
The space of metrics compatible with a given projective structures
can have dimensions $0, 1, 2, 3, 4$ or $6$. 
\end{theorem}
Our approach to the Koenigs's theorem is similar to that of Kruglikov's
\cite{Kruglikov}
who has however constructed an additional set of
invariants determining whether a metrisable projective structure
admits more than one metric in its projective class.

\section{Examples}
It is possible that the determinant (\ref{DEinvariant}) vanishes and
the projective structure $[\Gamma]$ is non metrisable either because
the further higher order obstructions do not vanish, or because
a solution to the Liouville system 
(\ref{linear_system})
is degenerate as a quadratic
form on $TU$. It can also happen when the 
projective structure fails to be real analytic.

In this section we shall give four examples illustrating 
this.
\label{examples}
\subsection{The importance of 6th order conditions}
Consider a one parameter family of 
homogeneous projective structures corresponding to 
the second order ODE
\[
\frac{d^2 y}{d x^2}=c\,  e^x+e^{-x}\Big(\frac{d y}{d x}\Big)^2.
\]
For generic $c$  the matrix ${\mathcal M}([\Gamma])$ has rank six and
the 5th order condition (\ref{DEinvariant1})  holds if $\hat{c}=48c-11$ 
is a root
of a quartic
\be
\label{quadric_5th}
\hat{c}^4-11286\, \hat{c}^2-850968\, \hat{c}-19529683=0.
\ee
The  6th order conditions (\ref{two_new_eq}) are satisfied iff
\[
3\,{\hat{c}}^{5}+ 529\,{\hat{c}}^{4}+222\,{\hat{c}}^{3}
-2131102\,{\hat{c}}^{2} -103196849\,{\hat{c}}-1977900451=0, 
\]
\[
\hat{c}^3-213\, \hat{c}^2-7849\, \hat{c}-19235=0.
\]

It is easy to verify that these three polynomials 
do not have a common root. Choosing $\hat{c}$ to be a real root
of (\ref{quadric_5th}) we can make the 5th order obstruction
(\ref{DEinvariant1}) vanish, but the two 6th order obstructions
$E_1, E_2$ do not vanish.
\subsection{The importance of the non-degenerate kernel}
This example illustrates why we cannot hope to characterise the
metrisability condition purely by vanishing of any set of invariants. 

Let $f$ be a smooth function on an open set $U\subset \R^2$.
Consider a one-parameter family of metrics 
\[
g_c=c \exp{(f(x,y))} dx^2 + dy^2,\qquad \mbox{where}\quad c\in \R^+. 
\]
The corresponding one-parameter family of projective structures
$[\Gamma_c]$  is given by the ODE
\[
\frac{d^2 y}{d x^2}
  = \frac{c}{2} \frac{\p f}{\p y} \exp{(f)}+ \frac{1}{2}\frac{\p  f}{\p x}
\Big(\frac{d y}{d x}\Big) 
+ \frac{\p f}{\p y} \Big(\frac{d y}{d x}\Big)^2.
\]
The 5th order obstruction (\ref{DEinvariant1}) 
and  6th order conditions $E_1, E_2$ of course vanish. 
Moreover rank ${\mathcal M}([\Gamma_c])=5$ for generic $f(x,y)$.

Now take the  limit $c=0$. The obstructions still vanish and 
rank ${\mathcal M}([\Gamma_0])=5$ but $[\Gamma_0]$ is not metrisable. 
This is because one can
select a 3 by 3 linear subsystem $\widetilde{\mathcal M}_0 \,\phi=0$, 
where $\phi=(\psi_1, \psi_2,
\mu)^T $,  from the 6 by 6 system
(\ref{geometric_form}). The 3 by 3 matrix $\widetilde{\mathcal M}_0$ can
be read off
(\ref{geometric_form}). For generic $f$
the  determinant of $\widetilde{\mathcal M}_0$ does not vanish and so there
does 
not exist a parallel  section ${\bf\Psi}$ of  (\ref{geometric_form}) such
that $\psi_1\psi_3-\psi_2^2\neq 0$. 
For example  $f=xy$ gives rank ${\mathcal M}([\Gamma_0])=5$  and 
\[
\det{(\widetilde{\mathcal M}_0)}=\frac{3xy}{4}-\frac{9}{2}.
\]

This non-metrisable example fails the genericity assumption 
$P([\Gamma])\neq0$
where $P([\Gamma])$ is given by (\ref{polynomial_adj}). The
kernel of  ${\mathcal M}([\Gamma_0])$  is spanned by a vector
$(0, 0, 1, 0, 0, 0)^T$ and the corresponding quadratic form
on $TU$ is degenerate.

\subsection{The importance of  real analyticity}
This example illustrates why we need to 
work in the real analytic  case to get sufficient conditions.  We shall
construct a simply connected projective surface in which every point has a 
neighbourhood on which there is a metric compatible with the given 
projective structure, but there is no metric defined on the whole 
surface that is compatible with the projective structure. 

Consider a plane $U=\R^2$ with cartesian coordinates $(x, y)$. Take two
constant coefficient metrics
on the  plane that are linearly independent, say, $g_+$ and $g_-$.  Now 
consider a modification of $g_-$ in the half-plane 
$x < -1$ such that the 
modified $g_-$ is the only global metric that is compatible with its 
underlying projective structure.  Similarly, consider a modification 
of $g_+$ on the half-plane $x > 1$ such that the modified $g_+$ is the only 
global metric that is compatible with its underlying projective structure.
The two projective structures agree (with the flat one) 
in the strip $-1 < x < 1$, so let the new projective structure be the 
one that agrees with that of modified $g_-$ 
when $x < 1$ and with the 
modified $g_+$ when $x > -1$.   This final projective structure will have 
compatible metrics locally near each point (sometimes, more than one, 
up to multiples), but will not have a compatible metric globally. Thus,
metrisability cannot be detected locally in the smooth category.

\subsection{One more  degenerate example.}
Take $\Gamma^2_{11}=A(x, y)$ and set all other components
of $\Gamma^a_{bc}$ to zero.
Equivalently, take $A_1=A_2=A_3=0, A_0=-A(x, y)$ (the case 
$A_0=A_1=A_2=0$ is also degenerate and can be obtained by reversing the role
of $x$ and $y$). 
For this degenerate case the Liouville relative invariant \cite{Liouville}
\begin{eqnarray*}
\nu_5&=&L_2(L_1\p_x L_2-L_2\p_xL_1)+L_1(L_2\p_y L_1-L_1\p_y L_2)\\
&&+A_3 (L_1)^3-A_2(L_1)^2 L_2+A_1 L_1(L_2)^2-A_0(L_2)^3
\end{eqnarray*}
vanishes.

The matrix ${\mathcal M}([\Gamma])$ in 
(\ref{DEinvariant})
has rank five and its determinant vanishes identically. In this case we can
nevertheless
analyse the linear system (\ref{linear_system}) directly without
even prolonging it. We solve for 
\[
\psi_2=-(1/2)y\alpha(x)+\beta(x)  ,\quad\psi_3=\alpha(x),
\]
where $\alpha$ and $\beta$ are some arbitrary functions of $x$,
and cross-differentiate the remaining equations to find
\begin{equation}
\label{intermedia_inv}
2\beta^{''}-y\alpha^{'''}+2(\p_x A)\alpha-2(\p_y A)\beta+(3A+y\p_y A)
\alpha'=0.
\end{equation}
Now assume further that $5\p_y^2 A+y\p_y^3 A\neq 0, \p_y^3 A\neq 0$
and perform further differentiations to eliminate
$\alpha, \beta$ from (\ref{intermedia_inv})
and to find the necessary metrisable condition for $A(x, y)$
\begin{eqnarray}
\label{degenerate1}
&&7(\p_y^3A)\;(\p_y^4A)\;(\p_x\p_y^3A)-5(\p_x\p_y^3A)\;(\p_y^5A)
\;(\p_y^2A)-
6(\p_x\p_y^4 A)\;(\p_y^3A)^2
\nonumber\\
&&+6(\p_y^5A)\;(\p_x\p_y^2A)\;(\p_y^3A)-7
(\p_y^4A)^2\;(\p_x\p_y^2A)+5(\p_x\p_y^4A)\;(\p_y^4A)\;(\p_y^2A)=0.
\end{eqnarray} 
The obstruction (\ref{degenerate1}) is of the same differential 
order as the 6 by  6 matrix (\ref{DEinvariant}), 
and we checked 
that it arises as a vanishing of a determinant of 
some  5 by 5 minors
of (\ref{DEinvariant}) (which factorise in this case with 
(\ref{degenerate1}) as a common factor).

We have pointed out that further genericity assumptions for $A$ were needed
to arrive at (\ref{degenerate1}). To construct an example of 
non-metrisable projective connection where these assumptions do not hold
consider the first Painlev\'e equation \cite{Ince}
\[
\frac{d^2 y}{d x^2}=6 y^2+x,
\]
for which  both (\ref{DEinvariant1}) and (\ref{degenerate1}) vanish.
However equation (\ref{intermedia_inv}) implies that
$\alpha(x)=\beta(x)=0$ so no metric exists in this case.
We would have reached the same conclusion by observing 
that in the  Painlev\'e I case rank$({\mathcal M})([\Gamma])=3$ and verifying
that the 6 by 6 matrix in the 8th order obstruction (\ref{rank_8}) has rank 5.
This obstruction therefore vanishes but the corresponding
one-dimensional kernel is spanned by
$(1, 0, 0, 0, 0, 0)^T$ and the corresponding solution to
the linear system (\ref{linear_system}) is degenerate.

In \cite{HD} it was shown that the Liouville invariant $\nu_5$
vanishes for all six Painlev\'e equations, and
we have verified that our invariant (\ref{DEinvariant1}) also
vanishes. The metrisability analysis would need to be done on
a case by case basis in a way analogous to our treatment
of Painlev\'e~I.
\section{Twistor Theory}
\label{sec_twistor}
In this Section we shall give  a twistorial treatment of the problem,
which clarifies the rather mysterious linearisation
(\ref{linear_system}) of the non-linear system 
(\ref{conditions}).

In the real analytic case one complexifies the projective structure, and
establishes a one-to-one correspondence between holomorphic 
projective structures $(U, [\Gamma])$ and
complex surfaces $Z$ with rational curves with self-intersection number 
one \cite{hitchin}. The points in $Z$ correspond to geodesics in $U$, and all
geodesics in $U$
passing through a point $u\in U$ form a rational curve $\hat{u}\subset Z$
with normal bundle $N(\hat{u})= {\mathcal O}(1)$. 
Here ${\mathcal O}(n)$ denotes the $n$th tensor power of the dual
of the tautological line bundle ${\mathcal O}(-1)$ over $\mathbb{P}(TU)$
which arises as a quotient of $TU-\{0\}$ by the Euler vector field.
Restricting the canonical line bundle $\kappa_Z$ of $Z$ 
to a twistor line $\hat{u}=\CP^1$ gives
\[
\kappa_Z=T^* (\hat{u})\otimes N^*(\hat{u})= {\mathcal O}(-3) 
\]
since the holomorphic tangent bundle to $\CP^1$ is 
${\mathcal O}(2)$. 
If
$U$ is a complex surface with a holomorphic projective structure, then its
twistor space $Z$ is $\mathbb{P}(TU)/D_x$, where
$D_x$ is the geodesic spray of the projective connection
(\ref{thomas_symbols})
\begin{eqnarray}
\label{spray_Dx}
D_x&=&z^a\frac{\partial}{\partial x^a} - 
\Pi_{ab}^c z^az^b\frac{\partial}{\partial z^c}\\
&=&\frac{\p}{\p x}+\zeta\frac{\p}{\p y}
+(A_0+\zeta A_1+\zeta^2A_2+\zeta^3 A_3)\frac{\p}{\p \zeta}.\nonumber
\end{eqnarray}
Here $(x^a, z^a)$ are coordinates on $TU$ and the second line uses projective
coordinate $\zeta=z^2/z^1$.
This leads to the double fibration
\[
U\longleftarrow \mathbb{P}(TU)\longrightarrow Z.
\] 
All these structures should be invariant under an anti-holomorphic
involution of $Z$ to recover a real structure on $U$. This works in the real
analytic case, but can in principle be extended to the smooth case
using the holomorphic discs of LeBrun-Mason \cite{LM}.

Now if the projective structure is metrisable,
$Z$ is equipped with a preferred section  of the anti-canonical divisor line
bundle ${\kappa_Z}^{-2/3}$ \cite{Cal,LM}.
The zero set of this section intersects each rational curve in $Z$ at two
points.
The pullback of this  section to $TU$ is a homogeneous function of
degree two $\sigma=\sigma_{ab}z^az^b$, where $z^a$ are homogeneous coordinates
on the fibres of $\mathbb{P}(TU)\rightarrow U$, and $\sigma_{ab}$ with 
$a, b=1, 2$ is a symmetric 2-tensor on $U$. 

This function Lie derives along the spray (\ref{spray_Dx})
and this gives the overdetermined linear system as the vanishing of a
polynomial homogeneous of degree 3 in $z^a$:
The condition $D_x(\sigma)=0$ implies the equation (\ref{twistor_eq})
which is equivalent to (\ref{linear_system}). 
\vspace{2ex}

We can understand the equation (\ref{twistor_eq})
using any connection $\Gamma$ in a projective class
instead of the projective connection $\nabla^{\Pi}$. To  see it we need 
to introduce a concept of projective weight \cite{Eastwood}.
First recall that the covariant derivative of the projective connection 
acting on vector fields is given by 
$\nabla_aX^c=\partial_aX^c+\Gamma_{ab}^cX^b$ and on $1$-forms by
$\nabla_a\phi_b=\partial_a\phi_b-\Gamma_{ab}^c\phi_c$. 
Let $\epsilon_{ab}=\epsilon_{[ab]}$ be a volume form on $U$. Changing
a representative of the projective class yields
\be
\label{vol_change}
\hat{\nabla}_a\epsilon_{bc}={\nabla}_a\epsilon_{bc}-3\;\omega_a\epsilon_{bc}.
\ee

Let $\mathcal{E}(1)$ be a line bundle over $U$ such that the
3rd power of its dual bundle is the canonical bundle of $U$.
The bundles  $\mathcal{E}(w)= \mathcal{E}(1)^{\otimes w}$
have a flat connection  induced from $[\Gamma]$. It changes
according to
\[
\hat{\nabla}_a h=\nabla_a h+ w\;\omega_a h
\]
under (\ref{equivalence}), where $h$ is a section of $\mathcal{E}(w)$.
\begin{definition}
\label{defi_weight}
The weighted vector field
with projective weight $w$  is a section of a bundle
$TU\otimes\mathcal{E}(w)$. 
\end{definition}
This definition naturally extends
to other tensor bundles. 
Now we shall choose a convenient normalisation
of $[\Gamma]$. For any choice of 
$\epsilon_{ab}$ we must have ${\nabla}_a\epsilon_{bc}=\theta_a\epsilon_{bc}$
for some $\theta_a$. We can change the projective representative
with $\omega_a=\theta_a/3$ 
and use
(\ref{vol_change}) to set $\theta_a=0$ so that $\epsilon_{ab}$ is parallel.
Let us assume that such a choice has been made.
We shall use the volume forms  to raise and lower indices according
to $z_a=\epsilon_{ba}z^b, z^a=z_b\epsilon^{ba}$. 
The residual freedom in (\ref{equivalence}) is to use $\omega_a=\nabla_a f
$ where $f$ is any function on $U$.  If $\nabla_a\epsilon_{bc}=0$ 
then \be
\label{vol_change2}
\hat{\nabla}_a\hat{\epsilon}_{bc}=0, \qquad \mbox{if}\quad 
\hat{\epsilon}_{ab}=e^{3f}\epsilon_{ab}.
\ee
Thus if $h\in \mathcal{E}(w)$ is a scalar of weight $w$ and we change the 
volume form as in (\ref{vol_change2}) then we must rescale
\[
h\longrightarrow \hat{h}=e^{wf} h
\]
with natural extension to other tensor bundles. Thus
$\epsilon^{ab}$ has weight $-3$.
\vspace{2ex}

Let us now come back to equation  (\ref{twistor_eq}) where the $\Pi$s 
are replaced by components of some connection in $[\Gamma]$ 
\[
\nabla_{(a}\sigma_{bc)}=0.
\]
If we change the representative of the projective class by 
(\ref{equivalence}) with $\omega_a=\nabla_a f$ the 
equation $D_x(\sigma)=0$ stays invariant if
\[
\sigma_{ab}\longrightarrow \hat{\sigma}_{ab}=e^{4f}\sigma_{ab}.
\]
This argument shows that
the linear operator 
\[
{\sigma}_{ab}\longrightarrow \nabla_{(a}\sigma_{bc)}
\]
is projectively invariant on symmetric two-tensors with weight 4.
Now $\sigma^{ab}:=\epsilon^{ac}\epsilon^{bd}\sigma_{cd}$ is a section
of $S^2(TU)\otimes \mathcal{E}(-2)$ and satisfies
\be
\label{tensor_prolong}
\nabla_a\sigma^{bc}=\delta_a^b\mu^c+\delta_a^c\mu^b
\ee
for some $\mu^b$.  The Liouville lemma \ref{Liouville_lemma} implies that 
if $\sigma^{ab}$ satisfies this equation then 
$g^{ab}=(\det{\sigma})\sigma^{ab}$ is a metric in the projective class.

The expression (\ref{tensor_prolong})
 is the tensor version of the first prolongation of the 
linear system (\ref{linear_system}). In the next section
we shall carry over the prolongation in the invariant manner
and express the 5th order obstruction (\ref{DEinvariant1})
as a weighted projective scalar invariant.
\section{An Alternative Derivation}
\label{sec_tractor}
In this section, we use the approach of \cite{East_Mat} to derive the
obstruction $\det{\mathcal{M}}([\Gamma])$ 
of Theorem~\ref{main_theorem}. One advantage of
this approach is that (\ref{DEinvariant1}) may then be 
written in terms of the
curvature of the connection and its covariant derivatives 
for any connection in
the given projective class. 
The symmetric form $\sigma^{ab}$ used in this Section 
is proportional to the quadratic form (\ref{metric})
and the objects $(\mu^a, \rho)$ are related but not equal to
$(\mu, \nu, \rho)$ defined by (\ref{munu}) and (\ref{Veq})
from Section \ref{sec_prolongation}.
Similarly the 6 by 6 matrix (\ref{sixbysix}) is related but not equal to
${\mathcal M}([\Gamma])$ given by (\ref{DEinvariant}). This is because
the choices made in the prolongation procedure leading 
(\ref{prolonged}) are different than those made in Section
\ref{sec_prolongation}. The resulting obstructions
(\ref{DEinvariant1}) and  (\ref{completecontraction})
do not depend on these choices and are the same up
to a non-zero exponential factor.

Let $\Gamma\in [\Gamma]$ be a connection in the projective class.
Its curvature is defined by
\[
[\nabla_a, \nabla_b]X^c=R_{abd}^c X^d
\]
and can be uniquely decomposed as
\be
\label{defofRho}
R_{abd}^c=
\delta_a^c\Rho_{bd}-\delta_b^c\Rho_{ad} +\beta_{ab}\delta^c_d
\ee
where $\beta_{ab}$  is skew. In dimensions higher than 2 there
would be another term (the Weyl tensor) in this curvature but
dimension in 2 it vanishes identically.

If we change the connection in the projective class
using (\ref{equivalence}) then
\[
\hat{\Rho}_{ab}={\Rho}_{ab}-\nabla_a\omega_b+\omega_a\omega_b,
\quad
\hat{\beta}_{ab}=\beta_{ab}+2\nabla_{[a}\omega_{b]}.
\]
The Bianchi identity implies that $\beta_{ab}$ is closed
and so locally it is clear that we can always choose a
connection in our projective class with $\beta_{ab}=0$
(in fact, this also true globally on an oriented manifold).
The residual freedom in changing the representative of the
equivalence class (1.1) is given by gradients
$\omega_a=\nabla_a f,$ where $f$ is a function on $U$.

Now  $\Rho_{ab}=\Rho_{ba}$ and the Ricci tensor of $\Gamma$ is symmetric.
The Bianchi identity  implies that  $\Gamma$ is flat
on a bundle of volume forms on $U$. Thus the  
normalisation of $\nabla_a$ may,
equivalently, be stated as requiring the existence of a volume form
$\epsilon^{ab}$ such that 
\[
\nabla_a\epsilon^{bc}=0.
\] 
Locally, such a volume
form is unique up to scale: let us fix one. 
This is the normalisation used in the previous Section.


The linear system and its prolongation
developed in \S\ref{sec_linearsystem} and \S\ref{sec_prolongation} is assembled
in \cite{East_Mat} into a single connection on a rank 
$6$ vector bundle over $U$.
Specifically, sections of this bundle comprise triples of contravariant tensors
$(\sigma^{ab},\mu^a,\rho)$ with $\sigma^{ab}$ being symmetric. The connection
is given by
\begin{equation}
\label{prolonged}\left\lgroup\begin{array}c
\sigma^{bc}\\[3pt]
\mu^b\\[3pt]
\rho
\end{array}\right\rgroup\stackrel{\nabla_a}{\longmapsto}
\left\lgroup\begin{array}c
\nabla_a\sigma^{bc}-\delta_a^b\mu^c-\delta_a^c\mu^b\\[3pt]
\nabla_a\mu^b-\delta_a^b\rho+\Rho_{ac}\sigma^{bc}\\[3pt]
\nabla_a\rho+2\Rho_{ab}\mu^b-2Y_{abc}\sigma^{bc}
\end{array}\right\rgroup,\end{equation}
where $Y_{abc}=\frac12(\nabla_a\Rho_{bc}-\nabla_b\Rho_{ac})$, the Cotton
tensor. The following is proved in~\cite{East_Mat}.
\begin{theorem} The connection $\nabla_a$ is projectively equivalent to a
Levi--Civita connection if and only if there is a covariantly constant section
$(\sigma^{ab},\mu^a,\rho)$ of the bundle with connection
{\rm(\ref{prolonged})} for which $\sigma^{ab}$ is non-degenerate.
\end{theorem}
It is also shown in \cite{East_Mat} how the rank $6$ bundle itself and its
connection (\ref{prolonged}) may be viewed as projectively invariant. In any
case, obstructions to the existence of a covariantly constant section may be
obtained from the curvature of this connection, which we now compute.
$$\begin{array}l\nabla_a\nabla_b\left\lgroup\begin{array}c
\sigma^{cd}\\[3pt]
\mu^c\\[3pt]
\rho
\end{array}\right\rgroup=
\nabla_a\left\lgroup\begin{array}c
\nabla_b\sigma^{cd}-\delta_b^c\mu^d-\delta_b^d\mu^c\\[3pt]
\nabla_b\mu^c-\delta_b^c\rho+\Rho_{bd}\sigma^{cd}\\[3pt]
\nabla_b\rho+2\Rho_{bc}\mu^c-2Y_{bcd}\sigma^{cd}
\end{array}\right\rgroup=\\[25pt]
\mbox{\small$\left\lgroup\begin{array}c
\nabla_a(\nabla_b\sigma^{cd}-\delta_b^c\mu^d-\delta_b^d\mu^c)
-\delta_a^c(\nabla_b\mu^d-\delta_b^d\rho+\Rho_{be}\sigma^{de})
-\delta_a^d(\nabla_b\mu^c-\delta_b^c\rho+\Rho_{be}\sigma^{ce})\\[3pt]
\nabla_a(\nabla_b\mu^c-\delta_b^c\rho+\Rho_{bd}\sigma^{cd})
-\delta_a^c(\nabla_b\rho+2\Rho_{bd}\mu^d-2Y_{bde}\sigma^{de})
+\Rho_{ad}(\nabla_b\sigma^{cd}-\delta_b^c\mu^d-\delta_b^d\mu^c)\\[3pt]
\nabla_a(\nabla_b\rho+2\Rho_{bc}\mu^c-2Y_{bcd}\sigma^{cd})
+2\Rho_{ac}(\nabla_b\mu^c-\delta_b^c\rho+\Rho_{bd}\sigma^{cd})
-2Y_{acd}(\nabla_b\sigma^{cd}-\delta_b^c\mu^d-\delta_b^d\mu^c)
\end{array}\right\rgroup$}\\[25pt]
\phantom{\nabla_a\nabla_b\left\lgroup\begin{array}c
\sigma^{cd}\\[3pt]
\mu^c\\[3pt]
\rho
\end{array}\right\rgroup}
=\left\lgroup\begin{array}c \nabla_a\nabla_b\sigma^{cd}
-\delta_a^c\Rho_{be}\sigma^{de}
-\delta_a^d\Rho_{be}\sigma^{ce}+\star\star\\[3pt]
\nabla_a\nabla_b\mu^c-\delta_a^c\Rho_{bd}\mu^d
+(\nabla_a\Rho_{bd})\sigma^{cd}+2\delta_a^cY_{bde}\sigma^{de}
+\star\star\\[3pt]
\nabla_a\nabla_b\rho+2(\nabla_a\Rho_{bc})\mu^c
-2(\nabla_aY_{bcd})\sigma^{cd}
+2Y_{abd}\mu^d+2Y_{acb}\mu^c+\star\star
\end{array}\right\rgroup,
\end{array}$$
where $\star\star$ denotes expressions that are manifestly symmetric in~$ab$.
Also notice that 
$$(\nabla_{[a}\Rho_{b]d})\sigma^{cd}+2\delta_{[a}^cY_{b]de}\sigma^{de}=
\delta_d^cY_{abe}\sigma^{de}+2\delta_{[a}^cY_{b]de}\sigma^{de}=
3\delta_{[a}^cY_{bd]e}\sigma^{de}=0,$$
and that 
$$Y_{[abc]}=0\implies Y_{acb}-Y_{bca}=Y_{abc}.$$
Therefore,
\begin{equation}\label{prolongedcurvature}
(\nabla_a\nabla_b-\nabla_b\nabla_a)\left\lgroup\begin{array}c
\sigma^{cd}\\[3pt]
\mu^c\\[3pt]
\rho
\end{array}\right\rgroup=
\left\lgroup\begin{array}c 0\\[3pt] 0\\[3pt]
10Y_{abc}\mu^c-4(\nabla_{[a}Y_{b]cd})\sigma^{cd}
\end{array}\right\rgroup.\end{equation}
Denoting the triple $(\sigma^{ab},\mu^b,\rho)$ by $\Sigma^\alpha$, we are
seeking a section $\Sigma^\alpha$ of our rank $6$ bundle so that
$\nabla_a\Sigma^\alpha=0$ and have found the explicit form of the evident
necessary condition $(\nabla_a\nabla_b-\nabla_b\nabla_a)\Sigma^\alpha=0$. 
We may rewrite our necessary
condition as $\epsilon^{ab}\nabla_a\nabla_b\Sigma^\alpha=0$. Notice, however,
that there is only one non-zero entry on the right hand side
of~(\ref{prolongedcurvature}). Our necessary condition analogous to
(\ref{linear_obstruction}) becomes
\be
\label{mikes_form}
\Xi_\alpha\Sigma^\alpha=0
\ee
for 
$$\Xi_\alpha\equiv\left\lgroup\begin{array}c
0\\[3pt]
5Y_a\\[3pt]
Z_{ab}
\end{array}\right\rgroup,\quad\mbox{where }
Y_c\equiv\epsilon^{ab}Y_{abc}\mbox{ and }
Z_{cd}\equiv-2\epsilon^{ab}\nabla_aY_{b(cd)}=\nabla_{(c}Y_{d)}.$$

Evidently, the quantity $\Xi_\alpha$ is a section of a rank $6$ bundle dual to
our previous one. Its sections consist of triples of covariant
tensors $(\kappa,\lambda_a,\tau_{ab})$ with $\tau_{ab}$ being symmetric and it 
inherits a connection dual to the previous one. Specifically,
$$\left\lgroup\begin{array}c
\kappa\\[3pt]
\lambda_b\\[3pt]
\tau_{bc}\\
\end{array}\right\rgroup\stackrel{\nabla_a}{\longmapsto}
\left\lgroup\begin{array}c
\nabla_a\kappa+\lambda_a\\[3pt]
\nabla_a\lambda_b+2\tau_{ab}-2\Rho_{ab}\kappa\\[3pt]
\nabla_a\tau_{bc}-\Rho_{a(b}\lambda_{c)}+2Y_{a(bc)}\kappa
\end{array}\right\rgroup,$$
where 
$$\left\lgroup\begin{array}c
\kappa\\[3pt]
\lambda_b\\[3pt]
\tau_{bc}\\
\end{array}\right\rgroup\intprod
\left\lgroup\begin{array}c
\sigma^{bc}\\[3pt]
\mu^b\\[3pt]
\rho
\end{array}\right\rgroup
\equiv\kappa\rho+\lambda_b\mu^b+\tau_{bc}\sigma^{bc}$$
is the dual pairing. By differentiating our necessary condition for
$\nabla_a\Sigma^\gamma=0$ we obtain
$$\Xi_\gamma\Sigma^\gamma=0\qquad(\nabla_a\Xi_\gamma)\Sigma^\gamma=0\qquad
(\nabla_{(a}\nabla_{b)}\Xi_\gamma)\Sigma^\gamma=0.$$
Since $\Sigma^\alpha$ is supposed to be a non-zero section, it follows that the
$6\times 6$ matrix 
\begin{equation}\label{sixbysix}\left\lgroup\left\lgroup\begin{array}c
0\\[3pt]
5Y_c\\[3pt]
Z_{cd}
\end{array}\right\rgroup\!\!,\:\nabla_a\!\!\left\lgroup\begin{array}c
0\\[3pt]
5Y_c\\[3pt]
Z_{cd}
\end{array}\right\rgroup\!\!,\:\nabla_{(a}\nabla_{b)}\!\!
\left\lgroup\begin{array}c
0\\[3pt]
5Y_c\\[3pt]
Z_{cd}
\end{array}\right\rgroup\right\rgroup\end{equation}
must be singular. Its determinant is the obstruction 
from Theorem~\ref{main_theorem}. We compute
$$\nabla_a\left\lgroup\begin{array}c
0\\[3pt]
5Y_c\\[3pt]
Z_{cd}\\
\end{array}\right\rgroup=
\left\lgroup\begin{array}c
5Y_a\\[3pt]
5\nabla_aY_c+2Z_{ac}\\[3pt]
\nabla_aZ_{cd}-5\Rho_{a(c}Y_{d)}
\end{array}\right\rgroup$$
and
$$\nabla_a\nabla_b\left\lgroup\begin{array}c
0\\[3pt]
5Y_c\\[3pt]
Z_{cd}\\
\end{array}\right\rgroup=
\left\lgroup\begin{array}c
5\nabla_aY_b+5\nabla_bY_a+2Z_{ba}\\[3pt]
\nabla_a(5\nabla_bY_c+2Z_{bc})+2\nabla_bZ_{ac}-10\Rho_{b(a}Y_{c)}
-10\Rho_{ac}Y_b\\[3pt]
\nabla_a(\nabla_bZ_{cd}-5\Rho_{b(c}Y_{d)})
-5\Rho_{a(c}\nabla_{|b|}Y_{d)}-2\Rho_{a(c}Z_{|b|d)}+10Y_{a(cd)}Y_b
\end{array}\right\rgroup$$
so
$$\nabla_{(a}\nabla_{b)}\left\lgroup\begin{array}c
0\\[3pt]
5Y_c\\[3pt]
Z_{cd}\\
\end{array}\right\rgroup=
\left\lgroup\begin{array}c
12Z_{ab}\\[3pt]
5\nabla_{(a}\nabla_{b)}Y_c+4\nabla_{(a}Z_{b)c}-5\Rho_{ab}Y_{c}
-15\Rho_{c(a}Y_{b)}\\[3pt]
\!\!\begin{array}r\nabla_{(a}\nabla_{b)}Z_{cd}-5(\nabla_{(a}\Rho_{b)(c})Y_{d)}
-5\Rho_{c(a}\nabla_{b)}Y_d-5\Rho_{d(a}\nabla_{b)}Y_c\quad\\
-\Rho_{c(a}Z_{b)d}-\Rho_{d(a}Z_{b)c}
+10Y_{(a}Y_{b)(cd)}\end{array}\!\!
\end{array}\right\rgroup.$$
To compute the determinant of the $6\times 6$ matrix (\ref{sixbysix}) we may
use the following.
\begin{lemma}\label{tedious} Let $\epsilon^{ab}$ denote the skew form in two
dimensions normalised as
$$\epsilon^{00}=0\quad\epsilon^{01}=1\quad
 \epsilon^{10}=-1\quad\epsilon^{11}=0.$$
Then the determinant of the $6\times 6$ matrix
$$\left\lgroup\begin{array}{cccccc}
     0&    P_0&    P_1&  Q_{00}&  Q_{01}&  Q_{11}\\
   R_0& S_{00}& S_{01}& T_{000}& T_{001}& T_{011}\\
   R_1& S_{10}& S_{11}& T_{100}& T_{101}& T_{111}\\
U_{00}&V_{000}&V_{001}&X_{0000}&X_{0001}&X_{0011}\\
U_{01}&V_{010}&V_{011}&X_{0100}&X_{0101}&X_{0111}\\
U_{11}&V_{110}&V_{111}&X_{1100}&X_{1101}&X_{1111}
\end{array}\right\rgroup$$
is
\begin{equation}\label{completecontraction}
\makebox[0pt]{$
\raisebox{30pt}{$\epsilon^{ab}\epsilon^{cd}\epsilon^{ef}\epsilon^{gh}
\epsilon^{ij}\epsilon^{kl}\epsilon^{mn}\epsilon^{pq}$}\!\!
\left[\!\!\!\begin{array}l
Q_{gi}S_{mp}T_{njk}U_{ac}V_{deq}X_{bfhl}-\frac16
P_pR_mS_{nq}X_{acgi}X_{behk}X_{dfjl}\\[3pt]
{}-\frac12
P_pS_{mq}T_{njl}U_{ce}X_{adgk}X_{bfhi}
-\frac12
P_pT_{mgi}T_{njk}U_{ac}V_{deq}X_{bfhl}\\[3pt]
\,{}+\frac12
P_pR_mT_{ngi}V_{acq}X_{dejk}X_{bfhl}
-\frac12
Q_{gi}R_mS_{np}V_{acq}X_{dejk}X_{bfhl}\\[3pt]
\,\,{}-\frac12
Q_{gi}R_mT_{njk}V_{acp}V_{deq}X_{bfhl}
-\frac14
Q_{gi}S_{mp}S_{nq}U_{ac}X_{dejk}X_{bfhl}\\[3pt]
\,\,\,{}-\frac14
Q_{gi}T_{mjk}T_{nhl}U_{ac}V_{dep}V_{bfq}
\end{array}\!\!\!\right]\enskip\quad$}\end{equation}
where $Q_{ab}=Q_{(ab)}$, $T_{cab}=T_{c(ab)}$, $U_{cd}=U_{(cd)}$,
$V_{cda}=V_{(cd)a}$, and $X_{cdab}=X_{(cd)(ab)}$.  
\end{lemma}
{\bf Proof.} A tedious computation.\koniec
Every tensor $Q, S, T, ...$ in this expression is constructed using 
one $\epsilon^{ab}$. Thus counting the total number
of $\epsilon^{ab}$s shows that the determinant
has a total projective weight $-42$
in a sense of Definition \ref{defi_weight}
which also means that it represents a section of the $14$th power of the
canonical bundle of $U$. 

To summarise, we have proved the following alternative formulation of
Theorem~\ref{main_theorem}.
\begin{theorem}
Suppose that $\nabla_a$ is a torsion-free connection in two-dimensions and that
$\epsilon^{bc}$ is a volume form such that $\nabla_a\epsilon^{bc}=0$. Define
the Schouten tensor\/ $\Rho_{ab}$ by {\rm(\ref{defofRho})} and
$$\textstyle Y_{abc}\equiv\frac12(\nabla_a\Rho_{bc}-\nabla_b\Rho_{ac})\qquad
Y_c\equiv\epsilon^{ab}\nabla_a\Rho_{bc}\qquad Z_{ab}\equiv\nabla_{(a}Y_{b)}.$$
Let
$$\begin{array}l
P_a\equiv 5Y_a\qquad 
Q_{ab}\equiv 12Z_{ab}\qquad
R_c\equiv 5Y_c\qquad
S_{ca}\equiv 5\nabla_aY_c+2Z_{ac}\\[3pt]
T_{cab}\equiv 5\nabla_{(a}\nabla_{b)}Y_c+4\nabla_{(a}Z_{b)c}-5\Rho_{ab}Y_{c}
-15\Rho_{c(a}Y_{b)}\\[3pt]
U_{cd}\equiv Z_{cd}\qquad
V_{cda}\equiv \nabla_aZ_{cd}-5\Rho_{a(c}Y_{d)}\\[3pt]
X_{cdab}\equiv \begin{array}[t]r
\nabla_{(a}\nabla_{b)}Z_{cd}-5(\nabla_{(a}\Rho_{b)(c})Y_{d)}
-5\Rho_{c(a}\nabla_{b)}Y_d-5\Rho_{d(a}\nabla_{b)}Y_c\\
-\Rho_{c(a}Z_{b)d}-\Rho_{d(a}Z_{b)c}
+10Y_{(a}Y_{b)(cd)}\end{array}
\end{array}$$
and define $\DG$ by the 
formula~{\rm(\ref{completecontraction})}. If $\nabla_a$ is projectively
equivalent to a Levi--Civita connection, then ${\DG}=0$.
\end{theorem}
In addition to giving an explicit formula for~$\DG$, there are
several other consequences of this theorem, which we shall now discuss. We have
found that 
$$\DG=\det \left\lgroup
\begin{array}{ccc}0&P_a&Q_{ab}\\
R_c&S_{ca}&T_{cab}\\
U_{cd}&V_{cda}&X_{cdab}
\end{array}\right\rgroup
$$
where
$$
Q_{ab}=Q_{(ab)}, \quad T_{cab}=T_{c(ab)}, \quad
U_{cd}=U_{(cd)},\quad V_{cda}=V_{(cd)a},\quad X_{cdab}=X_{(cd)(ab)}
$$
and the precise meaning of determinant is given by Lemma~\ref{tedious}. Though 
it makes no difference to the determinant and seemingly gives a more
complicated expression, it is more convenient to write 
$${\DG}=\frac{1}{4320}\det\bar{\Theta}\quad\mbox{where}\quad
\bar\Theta\equiv\left\lgroup\begin{array}{ccc}0&12P_a&Q_{ab}\\
30R_c&12S_{ca}&T_{cab}-5\Rho_{ab}R_c\\
30U_{cd}&12V_{cda}&X_{cdab}-5\Rho_{ab}U_{cd}
\end{array}\right\rgroup,$$
where the underlying matrix is evidently obtained by column operations from
the previous one. The reason is that this matrix better transforms under
projective change of connection. Specifically, if we write 
$$\bar\Theta=\left\lgroup\begin{array}{lll}0&\bar{P}_a&\bar{Q}_{ab}\\
\bar{R}_c&\bar{S}_{ca}&\bar{T}_{cab}\\
\bar{U}_{cd}&\bar{V}_{cda}&\bar{X}_{cdab}\end{array}\right\rgroup,$$
then under the change in connection 
$$\widehat\nabla_a\phi_b=\nabla_a\phi_b-\omega_a\phi_b-\omega_b\phi_a$$
induced by~(\ref{equivalence}), we find  
\begin{equation}\label{colops}\widehat{\bar\Theta}=
\left\lgroup\begin{array}{lll}
0&\tilde{P}_a&\tilde{Q}_{ab}-\tilde{P}_{(a}\omega_{b)}\\
\tilde{R}_c&\tilde{S}_{ca}-2\tilde{R}_c\omega_a&
\tilde{T}_{cab}-\tilde{S}_{c(a}\omega_{b)}+\tilde{R}_c\omega_a\omega_b\\
\tilde{U}_{cd}&
\tilde{V}_{cda}-2\tilde{U}_{cd}\omega_a&
\tilde{X}_{cdab}-\tilde{V}_{cd(a}\omega_{b)}+\tilde{U}_{cd}\omega_a\omega_b
\end{array}\right\rgroup\end{equation}
where 
\begin{equation}\label{rowops}
\begin{array}l\widetilde\Theta=\left\lgroup\begin{array}{lll}
0&\tilde{P}_a&\tilde{Q}_{ab}\\
\tilde{R}_c&\tilde{S}_{ca}&
\tilde{T}_{cab}\\
\tilde{U}_{cd}&
\tilde{V}_{cda}&
\tilde{X}_{cdab}
\end{array}\right\rgroup\\
=\left\lgroup\begin{array}{lll}0&\bar{P}_a&\bar{Q}_{ab}\\
\bar{R}_c&\bar{S}_{ca}-2\omega_c\bar{P}_a&
\bar{T}_{cab}-2\omega_c\bar{Q}_{ab}\\
\bar{U}_{cd}-\omega_{(c}\bar{R}_{d)}&
\bar{V}_{cda}-\omega_{(c}\bar{S}_{d)a}+\omega_c\omega_d\bar{P}_a&
\bar{X}_{cdab}-\omega_{(c}\bar{T}_{d)ab}+\omega_c\omega_c\bar{Q}_{ab}
\end{array}\right\rgroup.\end{array}\end{equation}
Notice that $\widetilde\Theta$ is obtained from $\bar\Theta$ by 
column operations and then $\widehat{\bar\Theta}$ is obtained from 
$\widetilde\Theta$ by row operations. It follows that determinant does 
not change, i.e.\ ${\mathcal{D}}(\widehat{\Gamma})=\DG$ is a projective 
invariant (from the formula (\ref{completecontraction}) it is already
apparent that $\DG$ is independent of choice of co\"ordinates). 
Thus we use the notation  ${\mathcal{D}}({[\Gamma]})$.

The argument following the formula (\ref{prolonged_operator}) shows that there
is only one 
obstruction to the metrisability at order 5 so  $\det{\mathcal{M}}([\Gamma])=0$ 
iff ${\mathcal{D}}({[\Gamma]})=0$. Thus
\[
\det{\mathcal{M}}([\Gamma])(d x\wedge dy)^{\otimes 14}
\]
is indeed a projective invariant as claimed in the Introduction.

A more invariant viewpoint on these matters is as follows. 
The formula (\ref{prolonged}) is for a connection on an invariantly defined
vector bundle,
denoted by ${\mathcal{E}}^{(BC)}$ in~\cite{East_Mat}. It arises from a
representation of ${\mathrm{SL}}(3,{\mathbb{R}})$ and the connection
(\ref{prolonged}) is closely related (but not equal to) the projective Cartan
connection induced on bundles so arising. The bundle is canonically filtered
with composition series
$${\mathcal{E}}^{AB}=
{\mathcal{E}}^{bc}(-2)+{\mathcal{E}}^b(-2)+{\mathcal{E}}(-2)$$
as detailed in~\cite{East_Mat}. Strictly speaking the quantity $\Xi_\gamma$ is 
not a section of the dual bundle ${\mathcal{E}}_{(CD)}$ but rather the
projectively weighted bundle ${\mathcal{E}}_{(CD)}(-5)$ with composition series
$${\mathcal{E}}_{(CD)}(-5)=
{\mathcal{E}}(-3)+{\mathcal{E}}_c(-3)+{\mathcal{E}}_{(cd)}(-3)$$
and $\bar\Theta$ is then obtained by applying the invariantly defined
`splitting operator' 
$${\mathcal{E}}(-5)\ni\xi\mapsto\left\lgroup\begin{array}c
30\xi\\
12\nabla_a\xi\\
\nabla_{(a}\nabla_{b)}\xi-5\Rho_{ab}\xi\end{array}\right\rgroup
\in{\mathcal{E}}_{(AB)}(-7)$$
coupled to the projectively invariant connection~(\ref{prolonged}). The upshot
is that $\bar\Theta$ is an invariantly defined section of
${\mathcal{E}}_{(CD)(AB)}(-7)$. Indeed, the formulae (\ref{colops}) and
(\ref{rowops}) giving $\widehat{\bar\Theta}$ in terms of $\bar\Theta$ are
precisely how sections of ${\mathcal{E}}_{(CD)(AB)}$ or
${\mathcal{E}}_{(CD)(AB)}(-7)$ transform under projective
change. Consequently, the obstruction ${\DG}$ is an invariant of
projective weight~$-42$.


\section{Outlook}
\label{sec_outlook}
In the  language of Cartan \cite{C22, BGH}, the general 2nd order
ODE (\ref{general_ODE}) defines a {\em path geometry}, and the paths are
geodesics
of projective connection if the ODE is of the form
(\ref{ODE2}). In this paper we have shown under what conditions
the paths in this geometry are unparametrised geodesics
of some metric. 
In case of higher dimensional projective structures the link with
ODEs is lost, but nevertheless one could search for 
conditions obstructing the metrisability in a way analogous to
what we did in Section (\ref{sec_tractor}). The results will have a different
character, however, owing to the presence of the Weyl curvature which will
modify the connection (\ref{prolonged}) as explained in \cite{East_Mat}.
The first necessary condition analogous to (\ref{DEinvariant1}) occurs already 
at order~2. Specifically, it is shown in \cite{East_Mat} that 
the curvature of the relevant connection in $n$ dimensions is given by 
$$(\nabla_a\nabla_b-\nabla_b\nabla_a)\left\lgroup\begin{array}c
\sigma^{cd}\\[3pt]
\mu^c\\[3pt]
\rho
\end{array}\right\rgroup=\left\lgroup\begin{array}c
W_{abe}^c\sigma^{de}+W_{abe}^d\sigma^{ce}
+\frac{2}{n}\delta_{[a}^cW_{b]ef}^d\sigma^{ef}
+\frac{2}{n}\delta_{[a}^dW_{b]ef}^c\sigma^{ef}
\\[3pt]
\ast\\[3pt]
\ast
\end{array}\right\rgroup$$
where $W_{abd}^c$ is the Weyl curvature and $\ast$ denotes  
expressions that we
shall not need. Since we are searching for covariant constant  
sections with
non-degenerate~$\sigma^{cd}$, in particular it follows that the linear
transformation $\sigma^{ef}\mapsto\Xi_{abef}^{cd}\sigma^{ef}$ where
\begin{equation}\label{defofXi}\textstyle\Xi_{abef}^{cd}:=
W_{ab(e}^c\delta_{f)}^d+
W_{ab(e}^d\delta_{f)}^c+
\frac{2}{n}\delta_{[a}^cW_{b](ef)}^d+
\frac{2}{n}\delta_{[a}^dW_{b](ef)}^c\end{equation}
is obliged to have a non-trivial kernel. Regarding $\Xi_{abef}^{cd}$  
as a
matrix representing this linear transformation, it should have $n(n 
+1)/2$
columns accounting for the symmetric indices $ef$. In its remaining  
indices it
is skew in $ab$, symmetric $cd$, and trace-free. These symmetries  
specify an
irreducible representation of ${\mathrm{GL}}(n,{\mathbb{R}})$ of  
dimension
$(n^2-1)(n^2-4)/4$, which we may regard as the number of rows of the  
matrix
$\Xi_{abef}^{cd}$. Notice that when $n=2$ this matrix is zero but as  
soon as
$n\geq3$ it has more rows than columns (for example, it is a $10 
\times 6$
matrix in dimension~3).
We claim that having a non-trivial kernel is a genuine condition and  
therefore
an obstruction to metrisability. For this, we need to show that
$\Xi_{abef}^{cd}$ can have maximal rank even when it is of the  
special form
(\ref{defofXi}) for some $W_{abd}^c$ having the symmetries of a Weyl  
tensor,
namely
\begin{equation}\label{symmsofW}
W_{abd}^c=-W_{bad}^c,\qquad W_{[abd]}^c=0,\qquad W_{abd}^a=0.
\end{equation}
Choose a frame and, for $n\geq 3$, consider the particular tensor $W_ 
{abd}^c$
having as its only non-zero components (no summation)
$$\begin{array}{ll}
W_{121}^1=-W_{211}^1=3(n^2-n-1)
& W_{122}^2=-W_{212}^2=3\\[3pt]
W_{123}^3=-W_{213}^3=-(n-1)(2n+3)
& W_{12c}^c=-W_{21c}^c=-(n-1),\,\forall c\geq 4\\[3pt]
W_{132}^3=-W_{312}^3=-(n^2-n-3)
& W_{1c2}^c=-W_{c12}^c=n+2,\,\forall c\geq 4\\[3pt]
W_{231}^3=-W_{321}^3=n(n+2)
& W_{2c1}^c=-W_{c21}^c=2n+1,\,\forall c\geq 4.
\end{array}$$
It is readily verified that the symmetries (\ref{symmsofW}) are  
satisfied. Form
the corresponding $\Xi_{abef}^{cd}$ according to (\ref{defofXi}) and  
consider
$\Xi_{12ef}^{cd}$. Being symmetric in $cd$ and~$ef$, we may regard it  
as a
square matrix of size $n(n+1)/2$ and it suffices to show that this  
matrix is
invertible. In fact, it is easy to check that it is diagonal with non- 
zero
entries along its diagonal.

\vspace{2ex}
The twistor analysis of Section (\ref{sec_twistor}) suggest
that there is some analogy between the metrisability problem we studied
in two dimensions
and existence of (possibly indefinite) K\"ahler
structure in a given anti-self-dual (ASD) conformal class ${\bf c}$ 
in on a four-manifold $M$.
A K\"ahler structure corresponds to a preferred section of
anti-canonical divisor
${\kappa_B}^{-1/2}$, where $\kappa_B$ is the canonical bundle of the 
twistor space \cite{pontecorvo} $B$ (a complex three-fold with
an embedded rational curve with normal bundle 
${\mathcal O}(1)\oplus{\mathcal O}(1)$).
Not all ASD structures are K\"ahler and the existence of the divisor
should lead to vanishing of some conformal invariants constructed out
of the ASD Weyl tensor (to the best of our knowledge they have never been 
written
down. If one adds the Ricci flat condition, some of the invariants are known
and can be expressed in terms of the Bach tensor).

These two constructions (ASD+K\"ahler in four dimensions
and projective + metrisable in two dimensions) are linked in the following way:
every ASD structure in $(2, 2)$ signature with a conformal null Killing
vector induces a projective structure on a two-dimensional space
$U$ of the $\beta$ surfaces (null ASD surfaces) in $M$. 
Conversely any two-dimensional projective structure
gives rise to (a class of) ASD structures with null conformal 
symmetry \cite{DW07, Cal}. Consider $B$ to be a holomorphic fibre bundle over
$Z$ with one dimensional fibres, where $Z$ is the twistor space of 
$(U,  [\Gamma])$ introduced in Section \ref{sec_twistor}.

Let $\hat{u}\subset Z$ be rational curve
in $Z$ corresponding to $u\in U$. The three-fold $B$ will be a twistor
space of an ASD conformal structure if $B$ restricts to ${\mathcal O}(1)$
on each twistor line $\hat{u}\subset Z$.
If $Z$ corresponds to a metrisable projective structure then the divisor
$\sigma$ lifts to a section of ${\kappa_B}^{-1/2}$, thus giving a $(2, 2)$ 
K\"ahler 
class. If the conformal Killing vector is not hyper-surface orthogonal
the local expression for the conformal class is
\[
{\bf c}=dz_a\otimes dx^a  - \Pi^c_{ab}\,z_c\, dx^a\otimes dx^b,
\]
where  $\Pi^c_{ab}$ are components of the projective connection
(\ref{thomas_symbols}).
The conformal Killing vector is a homothety $z_a/\p z_a$. 
This formula for ${\bf c}$  is equivalent to a special case
of expression (1.3) in \cite{DW07} 
after a change of coordinates and a conformal
rescaling (set $z_a=(-z e^t, e^t)$ and take $G=z^2/2+\gamma(x, y)z+\delta(x,y)$
for certain $\gamma, \delta$ in \cite{DW07}). It is a projectively invariant
modification of the Riemannian extensions of spaces with affine connection 
studied
by Walker \cite{Walker}.
The conformal class ${\bf c}$ is conformally flat
iff $[\Gamma]$ is projectively flat, i.e. its Cotton tensor
vanishes. This in turn is equivalent to the vanishing of the Liouville
expressions (\ref{liouville_inv}).

The  metrisable projective structures will therefore give rise
to (2, 2) ASD K\"ahler metric with conformal null symmetry.
Ultimately, the metrisability invariant (\ref{DEinvariant1}) in 
two dimensions will have its counterpart:  
a conformal invariant in four dimensions. Some progress in this direction 
has been made in \cite{DT09}.

\section*{Appendix}
\appendix
\def\theequation{\thesection{A}\arabic{equation}}

The connection $D=d+{\bfA}_1 dx+{\bfA}_2 dy$  
on the rank six vector bundle $\mathbb{E}\rightarrow U$ is
\begin{eqnarray}
\label{connection1}
{\bfA}_1&=&\left (
\begin{array}{cccccc}
-\frac{2}{3}A_1&2A_0 &0 &0 &0 &0\\
0&0 &0 &-\frac{1}{2} &0 &0\\
-2A_3&-\frac{2}{3}A_2 &\frac{4}{3}A_1 &0 &1 &0\\
(\bfA_1)_{41}&(\bfA_1)_{42} & (\bfA_1)_{43}& -\frac{1}{3}A_1& -3A_0&0\\
0&0 & 0& 0& 0&-1\\
(\bfA_1)_{61}&(\bfA_1)_{62}&(\bfA_1)_{63}&(\bfA_1)_{64} &(\bfA_1)_{65} 
&(\bfA_1)_{66}  
\end{array}
\right ),\nonumber \\
{\bfA}_2&=&
\left(
\begin{array}{cccccc}
-\frac{4}{3}A_2&\frac{2}{3}A_1 &2A_0 &1 &0 &0\\
0&0 &0 &0 &-\frac{1}{2}&0\\
0&-2A_3&\frac{2}{3}A_2 &0 &0 &0\\
0&0 &0 &0 &0 &-1\\
(\bfA_2)_{51}&(\bfA_2)_{52} &(\bfA_2)_{53} & 3A_3& \frac{1}{3}A_2&0\\
(\bfA_2)_{61}&(\bfA_2)_{62} &(\bfA_2)_{63} &(\bfA_2)_{64} &(\bfA_2)_{65} &
(\bfA_2)_{66}   
\end{array}
\right ). 
\end{eqnarray}
Let $V_1, ..., V_6$ be given by (\ref{formula_for_V}).
The curvature of $D$ is  
\begin{equation}
\label{six-curvature}
{\bf F}=d{\bfA}+{\bfA}\wedge {\bfA}=F\,d x\wedge d y=
\left(
\begin{array}{cccccc}
0&0&0&0&0&0\\
0&0&0&0&0&0\\
0&0&0&0&0&0\\
0&0&0&0&0&0\\
0&0&0&0&0&0\\
V_1&V_2&V_3&V_4&V_5&V_6
\end{array}
\right )\,d x\wedge d y. 
\end{equation}
The matrix elements of the connection are
\begin{eqnarray*}
(\bfA_1)_{41}&=&-\frac{4}{3}\p_xA_2+4A_0A_3+\frac{2}{3}\p_yA_1,\\
(\bfA_1)_{42}&=&-2\p_yA_0+\frac{2}{3}\p_xA_1+4A_2A_0-\frac{4}{9}(A_1)^2,\\
(\bfA_1)_{43}&=&2\p_xA_0-4A_0A_1,\\
(\bfA_1)_{61}&=& -\frac{4}{3}\p_x\p_y A_2-\frac{20}{3}A_0A_2A_3+\frac{2}{3}
(\p_y)^2A_1+4A_3\p_yA_0-2A_0\p_y A_3\\
&&-\frac{16}{9}A_2\p_x A_2+\frac{8}{9}A_2\p_y A_1  ,\\
(\bfA_1)_{62}&=& \frac{2}{3}\p_x\p_y A_1-\frac{4}{3}A_1\p_y A_1+2 A_0\p_y A_2-2 
(\p_y)^2 A_0+4 A_2\p_y A_0+4 A_3 \p_x A_0\\
&&+6A_0\p_x A_3+\frac{8}{9}A_1\p_x A_2
+\frac{4}{3}A_0A_1A_3-\frac{4}{3}A_0(A_2)^2 ,\\
(\bfA_1)_{63}&=& 2\p_x\p_y A_0+\frac{2}{3}A_0\p_x A_2-\frac{4}{3}A_2\p_x A_0
-4A_1\p_y A_0-\frac{4}{3}A_0\p_y A_1\\&&+\frac{8}{3}A_0A_1A_2+4 A_3 (A_0)^2,\\
(\bfA_1)_{64}&=&\frac{4}{3}\p_x A_2-\p_y A_1+5 A_0A_3,\\
(\bfA_1)_{65}&=&\frac{1}{3}\p_x A_1-4\p_y A_0
+3A_2A_0-\frac{2}{9}(A_1)^2,\\
(\bfA_1)_{66}&=&-\frac{1}{3}A_1,
\end{eqnarray*}
\begin{eqnarray*}
(\bfA_2)_{51}&=&-2\p_yA_3-4A_3A_2,\\
(\bfA_2)_{52}&=&2\p_xA_3-\frac{2}{3}\p_yA_2+4A_1A_3-\frac{4}{9}(A_2)^2,\\
(\bfA_2)_{53}&=&\frac{4}{3}\p_y A_1-\frac{2}{3}\p_x A_2+4A_0A_3,\\
(\bfA_2)_{61}&=&-2\p_x\p_y A_3-4A_2\p_xA_3-\frac{4}{3}A_3\p_x A_2
+\frac{2}{3}A_3\p_y A_1-\frac{4}{3}A_1\p_y A_3\\
&&-4A_0(A_3)^2-\frac{8}{3}A_1A_2A_3, \\
(\bfA_2)_{62}&=&2(\p_x)^2A_3-\frac{4}{3}A_2\p_x A_2-\frac{4}{3}A_0A_2A_3
-\frac{2}{3}\p_x\p_y A_2\\
&&+
4A_1\p_x A_3+4 A_0\p_y A_3+6 A_3\p_y A_0+\frac{4}{3}A_3(A_1)^2+2A_3\p_x A_1+
\frac{8}{9}A_2\p_y A_1,\\
(\bfA_2)_{63}&=&-\frac{2}{3}(\p_x)^2A_2+\frac{8}{9}A_1\p_x A_2
+4 A_0\p_x A_3-2 A_3\p_x A_0-\frac{16}{9}A_1\p_y A_1\\
&&+\frac{4}{3}\p_x\p_y A_1+\frac{20}{3}A_0A_1A_3,\\(\bfA_2)_{64}&=&
4\p_x A_3-\frac{1}{3}\p_y A_2+3 A_1A_3-\frac{2}{9}(A_2)^2,\\
(\bfA_2)_{65}&=& \p_x A_2-\frac{4}{3}\p_y A_1+5 A_0A_3, \\
(\bfA_2)_{66}&=&\frac{1}{3}A_2.\\
\end{eqnarray*}
The prolongation formulae are
\begin{eqnarray}
\label{PandQ}
P&=&-(\bfA_1)_{41}\,\psi_1-(\bfA_1)_{42}\,\psi_2-(\bfA_1)_{43}\,\psi_3
-(\bfA_1)_{44}\,\mu-
(\bfA_1)_{45}\,\nu,\nonumber\\
Q&=&-(\bfA_2)_{51}\,\psi_1-(\bfA_2)_{52}\,\psi_2-(\bfA_2)_{53}\,\psi_3
-(\bfA_2)_{54}\,\mu-
(\bfA_2)_{55}\,\nu,\nonumber\\
R&=&-(\bfA_1)_{61}\,\psi_1-(\bfA_1)_{62}\,\psi_2-(\bfA_1)_{63}\,\psi_3
-(\bfA_1)_{64}\,\mu-
(\bfA_1)_{65}\,\nu-(\bfA_1)_{66}\,\rho,\nonumber\\
S&=&-(\bfA_2)_{61}\,\psi_1-(\bfA_2)_{62}\,\psi_2-(\bfA_2)_{63}\,\psi_3
-(\bfA_2)_{64}\,\mu-
(\bfA_2)_{65}\,\nu-(\bfA_2)_{66}\,\rho.
\end{eqnarray}

\end{document}